


\newdimen\firstident
\newdimen\secondident
\newdimen\thirdident
\newbox\identbox
\newcount\assertioncount


\font\fractur=cmfrak


\catcode`\@=11




\let\pz=\S

%

%
%
%
%

\font\tenss=cmss10
\newfam\ssfam %
\textfont\ssfam=\tenss
\catcode`\_=11
\def\suf_fix{}
\def\scaled_rm_box#1{%
 \relax
 \ifmmode
   \mathchoice
    {\hbox{\tenrm #1}}%
    {\hbox{\tenrm #1}}%
    {\hbox{\sevenrm #1}}%
    {\hbox{\fiverm #1}}%
 \else
  \hbox{\tenrm #1}%
 \fi}
\def\suf_fix_def#1#2{\expandafter\def\csname#1\suf_fix\endcsname{#2}}
\def\I_Buchstabe#1#2#3{%
 \suf_fix_def{#1}{\scaled_rm_box{I\hskip-0.#2#3em #1}}
}
\def\rule_Buchstabe#1#2#3#4{%
 \suf_fix_def{#1}{%
  \scaled_rm_box{%
   \hbox{%
    #1%
    \hskip-0.#2em%
    \lower-0.#3ex\hbox{\vrule height1.#4ex width0.07em }%
   }%
   \hskip0.50em%
  }%
 }%
}
\I_Buchstabe B22
\rule_Buchstabe C51{34}
\I_Buchstabe D22
\I_Buchstabe E22
\I_Buchstabe F22
\rule_Buchstabe G{525}{081}4
\I_Buchstabe H22
\I_Buchstabe I20
\I_Buchstabe K22
\I_Buchstabe L20
\I_Buchstabe M{20em }{I\hskip-0.35}
\I_Buchstabe N{20em }{I\hskip-0.35}
\rule_Buchstabe O{525}{095}{45}
\I_Buchstabe P20
\rule_Buchstabe Q{525}{097}{47}
\I_Buchstabe R21 
\rule_Buchstabe U{45}{02}{54}
\suf_fix_def{Z}{\scaled_rm_box{Z\hskip-0.38em Z}}
\catcode`\"=12
\newcount\math_char_code
\def\suf_fix_math_chars_def#1{%
 \ifcat#1A
  \expandafter\math_char_code\expandafter=\suf_fix_fam
  \multiply\math_char_code by 256
  \advance\math_char_code by `#1
  \expandafter\mathchardef\csname#1\suf_fix\endcsname=\math_char_code
  \let\next=\suf_fix_math_chars_def
 \else
  \let\next=\relax
 \fi
 \next}
%
%
%
%
\def\font_fam_suf_fix#1#2 #3 {%
 \def\suf_fix{#2}
 \def\suf_fix_fam{#1}
 \suf_fix_math_chars_def #3.
}
\font_fam_suf_fix
 0rm
 ABCDEFGHIJKLMNOPQRSTUVWXYZabcdefghijklmnopqrstuvwxyz
\font_fam_suf_fix
 2scr
 ABCDEFGHIJKLMNOPQRSTUVWXYZ
\font_fam_suf_fix
 \slfam sl
 ABCDEFGHIJKLMNOPQRSTUVWXYZabcdefghijklmnopqrstuvwxyz
\font_fam_suf_fix
 \bffam bf
 ABCDEFGHIJKLMNOPQRSTUVWXYZabcdefghijklmnopqrstuvwxyz
\font_fam_suf_fix
 \ttfam tt
 ABCDEFGHIJKLMNOPQRSTUVWXYZabcdefghijklmnopqrstuvwxyz
\font_fam_suf_fix
 \ssfam
 ss
 ABCDEFGHIJKLMNOPQRSTUVWXYZabcdefgijklmnopqrstuwxyz
\catcode`\_=8
%
%
%
%
%
%
%
%
%
%
%
%
%
%
%
%
%
%
%
%
%
%
%
%
%
%

%

\def\hexnumber@#1{\ifcase#1 0\or1\or2\or3\or4\or5\or6\or7\or8\or9\or
        A\or B\or C\or D\or E\or F\fi }
\font\tenmsa=msam10  \font\tenmsb=msbm10
 
\font\sevenmsa=msam7 \font\sevenmsb=msbm7
   
\font\fivemsa=msam5  \font\fivemsb=msbm5
\newfam\msafam
\textfont\msafam=\tenmsa  \scriptfont\msafam=\sevenmsa
  \scriptscriptfont\msafam=\fivemsa
\edef\msa@{\hexnumber@\msafam}
\newfam\msbfam
\textfont\msbfam=\tenmsb  \scriptfont\msbfam=\sevenmsb
  \scriptscriptfont\msbfam=\fivemsb
\edef\msb@{\hexnumber@\msbfam}
\def\Bbb{\ifmmode\let\next\Bbb@\else
        \def\next{\errmessage{Use \string\Bbb\space only in math 
mode}}\fi\next}
\def\Bbb@#1{{\Bbb@@{#1}}}
\def\Bbb@@#1{\fam\msbfam#1}
\mathchardef\square="0\msa@03
\mathchardef\subsetneq="3\msb@28
\mathchardef\myleq="3\msa@36
\mathchardef\mygeq="3\msa@3E
\let\leq=\myleq
\let\geq=\mygeq


\catcode`\@=12

%

\def\C{{\Bbb C}}

\def\F{{\Bbb F}}
\def\G{{\Bbb G}}

\def\I{{\Bbb I}}

\def\P{{\Bbb P}}
\def\Q{{\Bbb Q}}
\def\R{{\Bbb R}}
\def\S{{\Bbb S}}

\def\Z{{\Bbb Z}}

%

\def\Abf{{\bf A}}

%

\def\hbar{{\bar h}}

\def\kbar{{\bar k}}

\def\sbar{{\bar s}}

%

\def\Ghat{{\hat G}}

%

%

%

\def\hline{\underline{h}}

\def\Mline{\underline{M}}

%

\def\Jtilde{{\tilde J}}

\def\Stilde{{\tilde S}}

\def\Xtilde{{\tilde X}}

\def\Ytilde{{\tilde Y}}

\def\Ztilde{{\tilde Z}}

%

%

%


%

%

\def\Fdbar{{\overline {\F}}}
\def\Qdbar{{\overline {\Q}}}

%

%

%

%

\def\hgbar{{\bar \eta}}

\def\lgbar{{\bar \lambda}}

%

%

%

%

%

%

%

%

%


\def\star{{}^*}

\def\cross{{}^{\times}}
\def\dlbrack{\lbrack{\mskip-2mu}\lbrack}      
\def\drbrack{\rbrack{\mskip-2mu}\rbrack}      
\def\lrangle{\langle\>{,}\>\rangle}           
\def\udot{{}^{\bullet}}
\def\Qp{\Q_{p}}
\def\Fp{{\F_{p}}}
\def\Qpbar{{{\overline {\Q}}_{p}}}

\def\Zp{\Z_{p}}


\def\overneq#1#2{\lower0.5pt\vbox{\lineskiplimit\maxdimen\lineskip-.5pt
                \ialign{$#1\hfil##\hfil$\crcr#2\crcr\not=\crcr}}}


\def\stimes{\mathbin{\raise1pt\hbox{$\scriptscriptstyle \bf
            \vert$}\mkern-5mu\times}}   


\def\ad{\mathord{\rm ad}}

\def\Aut{\mathord{\rm Aut}}

\def\codim{\mathord{\rm codim}}

\def\dim{\mathop{\rm dim}\nolimits}
\def\End{\mathop{\rm End}\nolimits}

\def\Frob{\mathop {\rm Frob}\nolimits}
\def\Gal{\mathop{\rm Gal}\nolimits}
\def\gr{\mathop{\rm gr}\nolimits}

\def\id{\mathop{\rm id}\nolimits}
\def\Im{\mathop{\rm Im}\nolimits}

\def\int{\mathop{\rm int}\nolimits}

\def\Ker{\mathop{\rm Ker}\nolimits}

\def\log{\mathop{\rm log}\nolimits}

\def\ord{\mathop{\rm ord}\nolimits}

\def\pr{\mathop{\rm pr}\nolimits}

\def\Res{\mathop{\rm Res}\nolimits}

\def\rk{\mathop{\rm rk}\nolimits}

\def\Sp{\mathop{\rm Sp}\nolimits}
\def\Span{\mathop{\rm Span}\nolimits}
\def\Spec{\mathop{\rm Spec}\nolimits}

\def\setback(#1){\mathrel{\mkern-#1mu}}
\def\varfill#1{$\smash{#1} \mkern-8mu
   \cleaders\hbox{$\mkern-3mu \smash{#1} \mkern-3mu$}\hfill
   \mkern-8mu #1$}
\def\equalfill{\varfill{=}}


\let\ar=\rightarrow

\let\alr=\leftrightarrow
\let\air=\hookrightarrow

\let\asr=\mapsto
\def\ariso{\buildrel\sim\over\ar}      

\let\arr=\longrightarrow
\let\all=\longleftarrow
\def\aerr{\arr\setback(27)\arr}        
\def\airr{\lhook\joinrel\arr}          

\def\arriso{\buildrel\sim\over\arr}    
\def\arrover#1{\buildrel#1\over\arr}   



\def\arrisoover#1{\mathop{\arr}\limits^{#1 \atop \sim}}


\def\arvar(#1){\hbox to #1pt{\rightarrowfill}}
\def\alvar(#1){\hbox to #1pt{\leftarrowfill}}
\def\arvarover(#1)#2{\mathop{\arvar(#1)}\limits^{#2}}
\def\alvarover(#1)#2{\mathop{\alvar(#1)}\limits^{#2}}


\let\ad=\downarrow

\def\add{\Big\ad}                      
\def\addleft#1{\llap{$\vcenter{\hbox{$\scriptstyle #1$}}$}\add}

\def\addright#1{\add\rlap{$\vcenter{\hbox{$\scriptstyle #1$}}$}}







\def\aqrvar(#1){\hbox to #1pt{\equalfill}}


\let\implies=\Rightarrow






%
%

%
%

%
%

\def\powerseries over #1 in #2{{#1 \dlbrack #2 \drbrack}}

%
%

\def\laurentseries over #1 in #2{{#1 (\!( #2 )\!)}}

%
%

\def\smallmatrix(#1,#2;#3,#4){\left({{#1\atop #3}\>{#2\atop #4}}\right)}

%
%

%
%

%
%

%
%

%

%

\def\restricted#1{\vert_{#1}}

%
%

\def\subscript#1\atop#2{_{\scriptstyle #1 \atop #2}}

%
%

%
%




\newfam\frac
\textfont\frac=\fractur
\scriptfont\frac=\fractur



\def\en_item#1#2{%
 \par
 \setbox\identbox=\hbox #1{#2}%
 \noindent
 \hangafter=1%
 \hangindent=\wd\identbox
 \box\identbox
 \ignorespaces
}

\def\ennopar_item#1#2{%
 \setbox\identbox=\hbox #1{#2}%
 \noindent
 \hangafter=1%
 \hangindent=\wd\identbox
 \box\identbox
 \ignorespaces
}

%
%

\def\hangitem#1{\en_item{}{#1\enspace}}

%
%
%

%

\def\hanghangitem#1{\en_item{}{\kern\firstident #1\enspace}}

%
%
\def\hanghanghangitem#1{\en_item{}{\kern\secondident #1\enspace}}


%
%
\def\indention#1{%
 \setbox\identbox =\hbox{\kern\parindent{#1}\enspace}%
 \firstident=\wd\identbox
}

%
%
\def\subindention#1{%
 \setbox\identbox=\hbox{{#1}\enspace}%
 \secondident=\wd\identbox
 \advance \secondident by \firstident
}

%
%
\def\subsubindention#1{%
 \setbox\identbox=\hbox{{#1\ }\enspace}%
 \thirdident=\wd\identbox
 \advance \thirdident by \secondident
}%

%
%
\def\litem#1{\en_item{to\firstident}{\kern\parindent#1\hfil \enspace }}
\def\llitem#1{\en_item{to\secondident}{\kern\firstident#1\hfil\enspace}}
\def\lllitem#1{\en_item{to\thirdident}{\kern\secondident#1\hfil\enspace}}
\def\ritem#1{\en_item{to\firstident}{\kern\parindent\hfil#1\enspace}}
\def\rritem#1{\en_item{to\secondident}{\kern\firstident\hfil#1\enspace}}
\def\rrritem#1{\en_item{to\thirdident}{\kern\secondident\hfil#1\enspace}}
\def\citem#1{\en_item{to\firstident}{\kern\parindent\hfil#1\hfil\enspace}}
\def\ccitem#1{\en_item{to\secondident}{\kern\firstident\hfil#1\hfil\enspace}}
\def\cccitem#1{\en_item{to\thirdident}{\kern\secondident\hfil#1\hfil\enspace}}
\def\rmlitem#1{\en_item{to\firstident}{\kern\parindent{\rm #1}\hfil \enspace }}
\def\rmllitem#1{\en_item{to\secondident}{\kern\firstident{\rm #1}\hfil\enspace}}
\def\rmlllitem#1{\en_item{to\thirdident}{\kern\secondident{\rm #1}\hfil\enspace}}

%
%

\def\lnoitem#1{\ennopar_item{to\firstident}{\kern\parindent#1\hfil \enspace }}
\def\llnoitem#1{\ennopar_item{to\secondident}{\kern\firstident#1\hfil\enspace}}
\def\lllnoitem#1{\ennopar_item{to\thirdident}{\kern\secondident#1\hfil\enspace}}

%
%

\def\assertionlist{
   \assertioncount=1
   \indention{(2)}}

\def\assertionitem{
   \litem{{\rm(\number\assertioncount)}}\advance\assertioncount by 1
}

%
%

\def\bulletlist{
   \indention{$\bullet$}}

\def\bulletitem{
   \litem{$\bullet$}}

%
%

%



\newcount\parno
\newcount\secno
\newcount\subsecno
\newdimen\partisize
\newdimen\parindentvalue


\hsize=146 true mm
\vsize=8.9 true in
\hoffset=0.6 true cm

\tolerance=300
\pretolerance=100

\tolerance=300
\pretolerance=100

\parindentvalue=15pt
\parindent=\parindentvalue

\mathsurround=0pt

\normallineskiplimit=.5pt
\normalbaselineskip=15pt
\normallineskip=1pt plus .5 pt minus .5 pt
\normalbaselines

\abovedisplayskip = 7pt plus 3pt minus 3pt
\abovedisplayshortskip = 1pt plus 2pt
\belowdisplayskip = 7pt plus 3pt minus 3pt
\belowdisplayshortskip = 3pt plus 2pt

\partisize=12pt


\global\parno=0
\global\secno=0
\global\subsecno=0


\newskip\parskipamount               
 \parskipamount=25pt plus4pt minus4pt
 \def\parskip{\vskip\parskipamount}
\newskip\secskipamount               
 \secskipamount=10pt plus4pt minus2pt
 \def\secskip{\vskip\secskipamount}


\font\partifont=cmbx7 at \partisize

\font\footnotefont=cmr8




\outer\def\paragraph#1{\goodbreak
   \global\secno=0
   \global\advance\parno by 1
   \parskip
   \hangitem{\partifont \the\parno~}
   {\partifont #1}
   \par\nobreak}

%

\outer\def\Notation#1: #2\par{\nobreak\secskip
  \itemitem{}{\bf Notation#1}:\enspace{\sl#2\par}\penalty50}

%

\outer\def\nnsection#1{\penalty-100
   \secskip\noindent
   {\bf #1}}

%

\outer\def\secstart#1{\penalty-100
   \global\subsecno=0
   \global\advance\secno by 1
   \secskip\noindent
   {\bf (\the\parno.\the\secno)~#1}}

%

\outer\def\proof#1: {\smallbreak {\it Proof\ #1}\/:\enspace}

%

\outer\def\claim#1: #2\par{\smallbreak
   {\bf #1}:\enspace{\sl#2\par}\penalty50}



%

\def\lformno{
   \global\advance\subsecno by 1
   \leqno(\the\parno.\the\secno.\the\subsecno)
   }


\def\displayno{
   \global\advance\subsecno by 1
   (\the\parno.\the\secno.\the\subsecno)
   }


\def\actualsecno{(\the\parno.\the\secno)}


\def\actualsubsecno{(\the\parno.\the\secno.\the\subsecno)}

%
%
%

\def\label#1{\xdef#1{\actualsecno}}

%
%
%

\def\sublabel#1{\xdef#1{\actualsubsecno}}



\outer\def\bibliography{
   \parskip
   \parindent=0pt
   {\partifont Bibliography}
   \vskip 18pt plus4pt minus2pt
   \frenchspacing
   }


\def\updot{{}^{\bullet}}
\def\dodot{{}_{\bullet}}
\def\upp{{}^{(p)}}
\def\DR{{\mathord{\rm dR}}}

\def\Ref{{\mathord{\rm Ref}}}
\def\Fzip{{\mathord{{\tt F}\hbox{{\rm -}}{\tt zip}}}}

\def\Fr{{\mathord{F}}}
\def\SP{{\mathord{\Sscr\Pscr}}}
\def\RP{{\mathord{\Rscr\Pscr}}}
\def\Par{{\mathord{\rm Par}}}
\def\Puni{{\mathord{\Pscr}}}
\def\Uuni{{\mathord{\Uscr}}}
\def\Fq{{\mathord{\F_q}}}
\def\OO{{\mathord{\rm O}}}
\def\relpos{{\mathord{\rm relpos}}}
\def\pr{{\mathord{\rm pr}}}
\def\sep{{\mathord{\rm sep}}}
\def\GL{{\mathord{\rm GL}}}
\def\SO{{\mathord{\rm SO}}}
\def\Stab{{\mathord{\rm Stab}}}
\let\phi=\varphi
\def\red{{\rm red}}

\def\allisounder#1{\mathop{\all}\limits^{\sim}_{#1}}
\def\FZ{{\mathord{\rm FZ}}}
\def\Fil{{\mathord{\rm Fil}}}
\catcode`\@=11
\mathchardef\nmid="3\msb@2D
\catcode`\@=12

\newcount\refno\refno=1
\def\incr{\advance\refno by 1}
\edef\BorTits{\number\refno}\incr
\edef\Bourb{\number\refno}\incr
\edef\DelCohInt{\number\refno}\incr
\edef\DelVarSh{\number\refno}\incr
\edef\DelRelKKK{\number\refno}\incr
\edef\DelIll{\number\refno}\incr
\edef\SGAthree{\number\refno}\incr
\edef\MAV{\number\refno}\incr
\edef\GeerKats{\number\refno}\incr
\edef\EGA{\number\refno}\incr
\edef\IllRay{\number\refno}\incr
\edef\KatoLog{\number\refno}\incr
\edef\KatzDiffEq{\number\refno}\incr
\edef\KottPts{\number\refno}\incr
\edef\Kraft{\number\refno}\incr
\edef\LusztPar{\number\refno}\incr
\edef\Mats{\number\refno}\incr
\edef\BMGSAS{\number\refno}\incr
\edef\BMSTPEL{\number\refno}\incr
\edef\BMDFEO{\number\refno}\incr
\edef\Oda{\number\refno}\incr
\edef\OgusKKK{\number\refno}\incr
\edef\OortStrat{\number\refno}\incr
\edef\OortTexel{\number\refno}\incr
\edef\RapZink{\number\refno}\incr
\edef\TWOrdin{\number\refno}\incr
\edef\TWDimStrat{\number\refno}\incr


\centerline{\partifont DISCRETE INVARIANTS OF VARIETIES}
\centerline{\partifont IN POSITIVE CHARACTERISTIC}
\vskip 3mm

\centerline{{\it by}}
\vskip 3mm

\centerline{Ben Moonen\footnote{*}{\footnotefont Research made possible by
a fellowship of the Royal Netherlands Academy of Arts and Sciences} and
Torsten Wedhorn}
\vskip 10mm

\noindent {\bf Abstract.} If $S$ is a scheme of characteristic $p$, we define
an $F$-zip over~$S$ to be a vector bundle with two filtrations plus a
collection of semi-linear isomorphisms between the graded pieces of the
filtrations. For every smooth proper morphism $X\to S$ satisfying certain
conditions the de Rham bundles $H^n_{{\rm dR}}(X/S)$ have a natural structure
of an $F$-zip. We give a complete classification of $F$-zips over
an algebraically closed field by studying a semi-linear variant of a variety
that appears in recent work of Lusztig. For every $F$-zip over $S$ our methods
give a scheme-theoretic stratification of $S$. If the $F$-zip is associated
to an abelian scheme over $S$ the underlying topological
stratification is the Ekedahl-Oort stratification. We conclude the paper with a
discussion of several examples such as good reductions of Shimura varieties of
PEL type and K3-surfaces.

\vskip 10mm

\noindent
{\partifont Introduction}

\secskip

Let $f \colon X \to S$ be a smooth proper morphism of schemes in
characteristic $p>0$. We say that $f$ satisfies condition (D) if the
sheaves $R^b f_* \Omega^a_{X/S}$ are locally free and if the Hodge-de
Rham spectral sequence degenerates at $E_1$-level. The de Rham
cohomology sheaves $M = H^m_{{\rm dR}}(X/S)$ are then locally free
${\cal O}_S$-modules that come naturally equipped with a descending
filtration~$C^\bullet$ (the Hodge filtration), an ascending
filtration~$D_\bullet$ (the conjugate filtration), and with ${\cal
O}_S$-linear isomorphisms $\phi_i\colon \bigl({\rm gr}^i_C\bigr)^{(p)}
\ariso {\rm gr}_i^D$ given by the (inverse) Cartier operator. We
call such a structure $\underline{M} =
(M,C^\bullet,D_\bullet,\phi_\bullet)$ an {\it $F$-zip\/} over~$S$.

In this paper we give a complete classification of $F$-zips over an
algebraically closed field. It turns out that $F$-zips over $k =
\bar{k}$ are essentially combinatorical objects. In order to state our
result, let us first define the {\it type\/} of an
$F$-zip~$\underline{M}$ over a connected basis~$S$ as the function
$\tau\colon {\Bbb Z} \to {\Bbb Z}_{\geq 0}$ given by $\tau(i) = {\rm
rank}_{{\cal O}_S} \bigl({\rm gr}^i_C\bigr)$. In the geometric example
considered above, the type is given by the Hodge numbers $h^{i,m-i}$
of the fibres of~$f$. Our first main result is the following
(cf.~section (4.4)).
\bigskip

\noindent
{\bf Theorem 1.} {\sl Let $k$ be an algebraically closed field of
characteristic $p>0$. Let $\tau \colon {\Bbb Z} \to {\Bbb Z}_{\geq 0}$
be a function with finite support $i_1 < \cdots < i_r$. Let $n_j :=
\tau(i_j)$, write $J = (n_r,\ldots,n_1)$, and let $n := n_1 + \cdots +
n_r$. Then there is a bijection
$$
\left\{\vcenter{
\setbox0=\hbox{isomorphism classes of $F$-zips}
\copy0
\hbox to\wd0{\hfil of type $\tau$ over $k$\hfil}}\right\}
\longleftrightarrow
(S_{n_r} \times \cdots \times S_{n_1})\backslash S_n =: {}^JW\, .
$$
More precisely, to each $u \in {}^JW$ we associate a ``standard
$F$-zip'' $\underline{M}^u_\tau$ over ${\Bbb F}_p$ such that any
$F$-zip~$\underline{M}$ over~$k$ is isomorphic to
$\underline{M}^u_\tau \otimes_{{\Bbb F}_p} k$ for some uniquely
determined $u \in {}^JW$.}
\bigskip

In the case of an abelian variety~$X$ over a perfect field~$k$, the
$F$-zip structure on $H^1_{{\rm dR}}(X/k)$ gives the Dieudonn\'e
module of the $p$-kernel group scheme $X[p]$. In this special case,
our classification theorem was proven (up to differences in
terminology) by Kraft in [14]. It was realized by Ekedahl and Oort
that this can be used to define a stratification of the moduli space
${\cal A}_g$ of abelian varieties in characteristic~$p$. This
Ekedahl-Oort stratification is a very useful tool in the study
of~${\cal A}_g$; see Oort, [22] and~[23].

Our theory of $F$-zips enables us to extend these ideas to arbitrary
families $f\colon X \to S$ satisfying condition~(D), and to de Rham
cohomology in arbitrary degree. We define a generalized Ekedahl-Oort
stratification of the base scheme~$S$. In fact, our theory gives a
natural scheme-theoretic definition of these strata, which is new even
in the case of abelian varieties. The result can be stated as
follows (cf.~section~(4.9)).
\bigskip

\noindent
{\bf Theorem 2.} {\sl Let $\tau$ and ${}^JW$ be as in Theorem~1. Let
$\underline{M} = (M,C^\bullet,D_\bullet,\phi_\bullet)$ be an $F$-zip
of type~$\tau$ over a scheme~$S$ of characteristic~$p$. For $u \in
{}^JW$ we define a subfunctor $S^u_{\underline{M}}$ of~$S$ by the
condition that a morphism $g\colon T \to S$ factors
through~$S^u_{\underline{M}}$ if and only if $g^*\underline{M}$ is
fppf-locally isomorphic to $\underline{M}^u_\tau \otimes_{{\Bbb F}_p}
{\cal O}_T$. Then $S^u_{\underline{M}} \subset S$ is representable by
a locally closed subscheme of~$S$, and the map
$$
\coprod_{u \in {}^JW} S^u_{\underline{M}}\; \airr S
$$
is a bijective monomorphism (i.e., a partition of~$S$).}
\bigskip

For the proof of Theorem~1 and Theorem~2 we study the
$\F_p$-scheme~$X_\tau$ whose $S$-valued points
are the triples $(C^\bullet,D_\bullet,\phi_\bullet)$ such
that $({\cal O}_S^n,C^\bullet,D_\bullet,\phi_\bullet)$ is an $F$-zip
of type~$\tau$ over~$S$. The algebraic group ${\rm GL}_{n,{\Bbb F}_p}$
naturally acts on~$X_\tau$. Theorem~1 amounts to a classification of
the ${\rm GL}_n$-orbits in~$X_\tau$.

We think of $X_\tau$ as a ``mod~$p$ analogue'' of a compactified
period domain. Indeed, if we let ${}^\#S \to S$ be the ${\rm
GL}_n$-torsor of trivialisations of the underlying vector bundle~$M$
then we get a natural ``mod~$p$ period map'' ${}^\#S \to X_\tau$,
analogous to the period maps arising in Hodge theory. It turns out
that there is a unique open ${\rm GL}_n$-orbit $X_\tau^{{\rm ord}}
\subset X_\tau$, the ``ordinary locus'' (cf.~section~(4.5)), which is to
be thought of as the interior of the period domain. In this picture,
the other strata correspond to degenerations of the data that
constitute an $F$-zip.
\medskip

In order to study the ${\rm GL}_n$-orbits in~$X_\tau$, we express the
latter in more group-theoretical terms. We introduce varieties~$Z_J$
that are semi-linear variants of the varieties studied by Lusztig
in~[15]. We consider these varieties in the general context of a (not
necessarily connected) reductive group~$G$ over a finite field. Write
$(W,I)$ for the Weyl group of~$G$ with its set of simple
reflections. As further input for the definition of~$Z_J$ we need two
subsets $J$, $K \subseteq I$, and a Weyl group element $x \in W$
satisfying certain assumptions (see section~(3.2)). Write~$U_P$ for
the unipotent radical of a parabolic $P \subset G$. Then $Z_J$ is the Zariski
sheafification of the functor that classifies triples
$\bigl(P,Q,[g]\bigr)$ with $P$ and~$Q$ parabolic subgroups of types
$J$ and~$K$, respectively, and with $[g]$ a double coset in
$U_Q\backslash G/F(U_P)$ such that $Q$ and~${}^gF(P)$ are in relative
position~$x$. We prove that $Z_J$ is a smooth variety of dimension
equal to $\dim(G)$. The group~$G$ naturally acts on~$Z_J$.

The connection with the theory of $F$-zips is as follows. Let $\tau$
and~$J$ be as in Theorem~1, and take $G = {\rm GL}_{n,{\Bbb F}_p}$. We
identify $W = S_n$. The ordered partition~$J = (n_r,\ldots,n_1)$
corresponds to a subset of the set~$I$ of simple reflections. For $K
\subseteq I$ we take the subset corresponding to the opposite
partition $(n_1,\ldots,n_r)$, and for $x$ we take the element of
minimal length in the double coset $W_K w_0 W_J$, where $W_J$ and $W_K
\subset W$ are the subgroups generated by~$J$ and~$K$, respectively,
and where $w_0 \in W$ is the longest element. We show that with these
choices, there is a ${\rm GL}_n$-equivariant isomorphism between
$X_\tau$ and the variety~$Z_J$. Theorem~1 is then a consequence of the
following general result about the varieties~$Z_J$ (cf.~section~(3.25)).
\bigskip

\noindent
{\bf Theorem 3.} {\sl There is a bijection between the set of
$G$-orbits in~$Z_J$ and the set ${}^JW \subset W$ of elements $w\in W$
that are of minimal length in their coset~$W_J w$. (So ${}^JW$ is in
bijection with $W_J\backslash W$.)}
\bigskip

The idea for the proof of this theorem is the following. Let
$\bigl(P,Q,[g]\bigr)$ be a point of~$Z_J$. We define a new pair of
parabolics $(P_1,Q_1)$ by
$$
P_1 := (P\cap Q)U_P\, ,\qquad Q_1 := (Q \cap gF(P_1)g^{-1})U_Q\, .
$$
In a sense that can be made precise, the pair $(P_1,Q_1)$ is a
refinement of the pair $(P,Q)$. Repeating this process, we get a
sequence of pairs $(P_n,Q_n)$ that stabilizes. Then the bijection in
Theorem~3 is obtained by sending the point $\bigl(P,Q,[g]\bigr)$ to
the element of~$W$ that measures the relative position of $P_n$
and~$Q_n$ for $n \gg 0$.
\medskip

The same ideas as sketched here can be applied to study $F$-zips with
certain additional structures, such as a bilinear form or an action of
a semi-simple algebra. We apply this to abelian varieties,
K3-surfaces, and to good reductions of PEL-Shimura varieties. In this
last case, we give a new proof of the dimension formula for
Ekedahl-Oort strata that was obtained in~[\BMDFEO] using the results
of~[\TWDimStrat]. In fact, this is a consequence of the following
general result on the dimensions of the $G$-orbits in~$Z_J$
(cf.~section~(3.20)).
\bigskip

\noindent
{\bf Theorem 4.} {\sl For $u \in {}^JW$, let $O^u \subset Z_J$ be the
corresponding $G$-orbit under the bijection of Theorem~3. Then
$$
{\rm codim}(O^u,Z_J) = \dim({\rm Par}_J) - \ell(u)\, ,
$$
where $\ell(u)$ is the length of~$u$ in the Coxeter group~$W$, and
where ${\rm Par}_J$ is the variety of parabolics of type~$J$.}
\medskip

We will now give an overview of the structure of the paper. In the
first section we give the definition of $F$-zips over a base scheme of
characteristic~$p$, and we define standard $F$-zips. Section~2
contains some notations and lemmas on parabolic subgroups of reductive
groups, their relative position, and their Levi subgroups.

Section~3 is the technical heart of the paper. Here we introduce and
study the varieties~$Z_J$ discussed above. The main goal of this
section is the study of the fppf-quotient $G\backslash Z_J$ and, as an
application, the proof of Theorems 3 and~4. Our method is a variation
on ideas of Lusztig in~[\LusztPar].

In Section~4 we prove Theorems~1 and~2 announced above. The proof is
an easy application of our classification of the $G$-orbits in~$Z_J$.

In Section~5 we briefly discuss some examples of $F$-zips with
additional structure. Finally, in Section~6 we discuss applications to
geometry. We explain how a morphism $X \to S$ satisfying condition~(D)
gives rise to an $F$-zip structure on the de Rham cohomology. There is
also a version of this for log-schemes. We show that it is possible to
detect ordinariness, in the sense of Illusie and Raynaud~[10], from our
partition of~$S$. Next we consider good reductions of Shimura
varieties of PEL-type. The partition obtained in this case is the
generalized Ekedahl-Oort stratification studied earlier by the authors
in [17], [18] and~[27]. Finally, we study K3-surfaces $X \to S$, and
we make the connection between the stratification of~$S$ given by the
height and the Artin invariant, and the generalized Ekedahl-Oort
stratification obtained by our methods.
\bigskip

\noindent
{\bf Acknowledgements.} The first author thanks the Mathematics
Institute of the University of Cologne for its support and hospitality
during a short visit in May 2002, when part of the research in this
paper was conducted.

The second author is grateful to the European Algebraic Geometry
Research Training Network (EAGER) for the financial support of a visit
in Amsterdam in November 2001 when the research for this paper was
started. He further thanks the University of Amsterdam for its support
and hospitality during two visits in November 2001 and April 2003.



\paragraph{Filtrations and Flags}

\secstart{}\label{\Frdef} Throughout this section, $p$ is a prime
number and $q$ is a fixed power of~$p$. For a scheme~$S$ of
characteristic~$p$ we denote by $\Fr_S\colon S \to S$ the morphism
which is the identity on the underlying topological space and the
homomorphism $x \asr x^q$ on the sheaves of rings. For an
$\Oscr_S$-module $M$ we set $M^{(q)} = \Fr_S^*M$.

\secstart{}\label{\defineFil} Let $S$ be a scheme, and let $M$ be
a locally free $\Oscr_S$-module of finite rank. By a {\it descending
filtration} $C\updot$ of~$M$ we mean a sequence $(C^i)_{i\in\Z}$ of
$\Oscr_S$-submodules $C^i \subset M$ such that $C^i$ is locally
on~$S$ a direct summand of~$C^{i-1}$ and such that $C^i = M$ for
$i \ll 0$ and $C^i = (0)$ for $i \gg 0$. We set $\gr^i_C(M) =
\gr^i_C = C^i/C^{i+1}$. We have an analogous definition of an
{\it ascending filtration}~$D\dodot$, with associated graded modules
$\gr_i^D(M) = \gr_i^D = D_i/D_{i-1}$.

\secstart{}\label{\defineFlag} Let $M$ be as above. A {\it flag\/}
of~$M$ is a set~$\Delta$ of $\Oscr_S$-submodules of~$M$ which are
locally direct summands, such that $\Delta$ contains $(0)$ and~$M$
and is totally ordered by inclusion. Every (descending or ascending)
filtration~$C\updot$ defines a flag by forgetting the enumeration.

The set of flags of~$M$ is partially ordered by inclusion. We say
that $\Delta$ is a {\it refinement\/} of~$\Delta'$ if $\Delta
\supset \Delta'$.

\secstart{}\label{\typelocconst} Let $S$ be a scheme and let
$C\updot$ be a descending filtration of a locally free
$\Oscr_S$-module~$M$ of finite type. For $s \in S$, consider the
function $\tau^s_{C\updot}\colon\Z \to \Z_{\geq 0}$ given by $m
\asr \dim_{\kappa(s)}(\gr^m_{C \otimes \kappa(s)})$. As the
$\gr^m_C$ are locally free, the function $\tau\colon s \asr
\tau^s_{C\updot}$ is locally constant; it takes values in the set
of maps $\Z \to \Z_{\geq 0}$ with finite support. We refer
to~$\tau$ as the {\it type of the filtration\/}. If $S$ is
connected then $\tau$ is given by a single function $\Z \to
\Z_{\geq 0}$ with finite support.

A similar definition applies to ascending filtrations.

\secstart{Definition}:\label{\defineFzip} Let $S$ be an
$\F_q$-scheme. An {\it $F$-zip over~$S$\/} is a tuple $\Mline =
(M,C\updot,D\dodot, \varphi\dodot)$ where
\indention{--}
\litem{--}
$M$ is a locally free $\Oscr_S$-module of finite rank,
\litem{--}
$C\updot$ is a descending filtration of~$M$,
\litem{--} 
$D\dodot$ is an ascending filtration of~$M$,
\litem{--}
$\varphi\dodot$ is a family of $\Oscr_S$-linear isomorphisms
$$\varphi_n\colon (\gr^n_C)^{(q)} \arriso \gr_n^D$$
for $n \in \Z$.

The rank of $M$ is called the {\it height\/} of~$\Mline$. The type
of the filtration~$C\updot$ is called the {\it type\/}
of~$\Mline$.

We have the obvious notion of a morphism of $F$-zips (morphisms are
not required to be strict for the filtrations) and hence get
the category of $F$-zips over~$S$, which is an $\F_q$-linear rigid
tensor category. Note that this category is not abelian.

\secstart{} For an $\Fq$-scheme~$S$, let $\Fzip(S)$ be the category
which has as objects the $F$-zips over~$S$ and as morphisms the
isomorphisms of $F$-zips over~$S$. For a morphism of schemes
$f\colon T \to S$ we have an obvious pullback functor $f^*\colon
\Fzip(S) \to \Fzip(T)$. In this way we obtain a stack $\Fzip$,
fibered over the category of $\F_q$-schemes endowed with the fpqc
topology.

\secstart{Proposition}: {\sl The stack $\Fzip$ is a smooth Artin
stack over~$\F_q$. If $\tau\colon \Z \to \Z_{\geq 0}$ is a
function with finite support, the substack $\Fzip^{\tau}$ of
$F$-zips of type~$\tau$ is an open and closed substack of~$\Fzip$
and we obtain a decomposition
$$
\Fzip = \coprod_{\tau} \Fzip^{\tau}\, .
$$
The Artin stacks $\Fzip^{\tau}$ are quasicompact.}
\smallbreak

The easy proof is omitted.

\secstart{Example}: Assume $q = p$ and let $S = \Spec(R)$ with $R$
a perfect ring of characteristic~$p$. Consider a
BT$_1$-Dieudonn\'e module over~$S$, by which we mean a triple
$(M,F,V)$ with $M$ a projective $R$-module of finite type,
$F\colon M \to M$ an $\Frob_R$-linear map, $V\colon M \to M$ an
$\Frob_R^{-1}$-linear map, such that $\Ker(F) = \Im(V)$ and $\Im(F)
= \Ker(V)$ are locally direct summands of~$M$.

The category of BT$_1$-Dieudonn\'e modules can be identified with
the category of $F$-zips whose type~$\tau$ has support
contained in $\{0,1\}$: To $(M,F,V)$ we associate the
$F$-zip $(M,C\updot,D\dodot,\varphi\dodot)$ with
$$
\displaylines{C^0 = M \;\supset\; C^1 = \Ker(F) \;\supset\; C_2 = (0) \cr
D_{-1} = (0) \;\subset\; D_0 = \Im(F)\ \;\subset\; D_1 = M\, , \cr}
$$
with $\varphi_0\colon \bigl(M/\Ker(F)\bigr)^{(p)} \arriso M$ the
(linearization of the) isomorphism induced by~$F$, and
$\varphi_1\colon \Ker(F)^{(p)} \arriso M/\Im(F)$ the inverse of
the (linearized) isomorphism induced by~$V$.

\secstart{Standard F-zips}:\label{\StandFzips} We fix an integer
$n \geq 1$ and a map $\tau\colon \Z \to \Z_{\geq 0}$ with
$\sum_{i\in \Z} \tau(i) = n$. Let $i_1 < \cdots < i_r$ be the
support of~$\tau$ and $J = (n_r,\dots,n_1)$ be the ordered
partition of~$n$ with $n_j = \tau(i_j)$. (Note the order of
the~$n_j$.) Let~$W = S_n$ be the group of permutations of
$\{1,\dots,n\}$ and consider $W_J = S_{n_r} \times \cdots \times
S_{n_1}$ as a subgroup of~$S_n$ in the usual way. We set $m_j =
n_1 + \cdots + n_j$. Let $x \in W$ be defined by
$$
x(i) = i + m_j + m_{j-1} - n\, ,\qquad\hbox{if $n - m_j < i \leq n
- m_{j-1}$}\, ,
$$
i.e., $x$ is the element of minimal length in $w_0W_J$, where
$w_0$ is the longest element in~$W$. Finally let ${}^JW$ be the
set of permutations $u \in W$ with the property that
$$
u^{-1}(n - m_j + 1) < u^{-1}(n - m_j + 2) < \cdots < u^{-1}(n - m_{j-1})
$$
for all $j = 1,\dots,r$, i.e., ${}^JW$ consists of those $u \in W$
which are of minimal length in their left coset~$W_Ju$.

To $\tau$ and $u \in {}^JW$ we associate a standard $F$-zip
$\Mline_\tau^u = \bigl(M^u_{\tau}, (C\updot)^u_{\tau},
(D\dodot)^u_{\tau}, (\varphi\dodot)^u_{\tau}\bigr)$ over~$\F_p$,
where \bulletlist \bulletitem $M^u_{\tau} = \F_p^n\, ,$
\bulletitem $(C\updot)^u_{\tau}$ is the unique filtration of
type~$\tau$ such that the associated flag is given by
$$
\eqalign{
\F_p^n \supset \F_p^{\{u(1),u(2),\dots,u(m_{r-1})\}} &\supset
\F_p^{\{u(1),u(2),\dots,u(m_{r-2})\}} \cr
&\supset \cdots \supset \F_p^{\{u(1),u(2),\dots,u(m_1)\}} \supset
(0)\, ,\cr
}
$$
\bulletitem $(D\dodot)^u_{\tau}$ is the unique filtration of type
$\tau$ such that the associated flag is given by
$$
(0) \subset \F_p^{\{1,\dots,n-m_{r-1}\}} \subset
\F_p^{\{1,\dots,n-m_{r-2}\}}
\subset \cdots \subset \F_p^{\{1,\dots,n-m_0\}} = \F_p^n\, ,
$$
\bulletitem $(\varphi_i)^u_{\tau}$ is zero for $i \notin
\{i_1,\dots,i_r\}$ and for $i = i_j$ it is the isomorphism
$$
(\gr^i_C)^{(p)} = \F_p^{\{u(m_{r-j} +1), \dots, u(m_{r-j+1})\}} \arriso
\gr_i^D = \F_p^{\{n-m_{r-j+1} + 1, \dots, n-m_{r-j}\}}
$$
induced by the permutation matrix associated to $x^{-1}u^{-1}$.



\paragraph{The relative position of parabolics over an arbitrary base}

\nnsection{}In this section we introduce some notations and
collect some facts about parabolics of a reductive group~$G$ over
an arbitrary base. At the end we explain all these notions for the
case $G = \GL_n$.

\secstart{} Let $G$ be a group, $X \subset G$ a subset and $g \in
G$. Then we set ${}^gX = gXg^{-1}$.

If $P$ is a parabolic subgroup of some reductive group scheme, we
denote by $U_P$ its unipotent radical.

\secstart{}\label{\loccstsch} Let $k$ be a field and let $k^\sep$ be a
separable closure of~$k$. Recall that the functor $X \asr X(k^\sep)$ gives
an equivalence of the category of (finite) \'etale $k$-schemes with the
category of finites discrete sets endowed with a continuous action of
$\Gal(k^\sep/k)$.

\secstart{}\label{\groupnotation} We fix the following notations:
Let $k$ be a field and let $k^\sep$ be a separable closure of~$k$.
We denote by~$S$ an arbitrary $k$-scheme. If $X$ is a $k$-scheme,
write $X_S := X \times_{k} S$.

Let $G$ be a connected quasi-split reductive group over~$k$. We
denote by $(W,I)$ the Weyl group of the abstract based root datum
of~$G$, together with its set of simple reflections. It is a
finite Coxeter system carrying a continuous action of $\Gal(k^\sep/k)$.

For subsets $J$, $K \subset I$ we denote by $W_J$ the subgroup
of~$W$ generated by~$J$ and by ${}^JW^K$ the set of elements $w
\in W$ that are of minimal length in their double coset $W_JwW_K$.
We write ${}^JW = {}^JW^{\emptyset}$ and $W^K =
{}^{\emptyset}W^K$.

\secstart{}\label{\parabolics} Let $\Par$ be the smooth proper
$k$-scheme that parametrizes the parabolic subgroups of~$G$. It
carries a $G$-action and the fppf quotient $G\backslash \Par$ is
representable by a finite \'etale $k$-scheme~$\Dscr$, see SGA3 ([\SGAthree],
Exp.~XXVI, section~3, where $\Dscr$ is called $\Pscr\bigl({\bf
Dyn}(G)\bigr)$). We denote by
$$
\tbf\colon \Par \arr \Dscr
$$
the canonical morphism. For $P \in \Par(S)$ we call $\tbf(P) \in
\Dscr(S)$ the {\it type of\/}~$P$. If $J$ is a section of~$\Dscr$
over~$S$ we denote by $\Par_J = \tbf^{-1}(J) \subset \Par_S$ the
scheme of parabolics of type~$J$.

Under the equivalence of~\loccstsch, the scheme $\Dscr$
corresponds to the powerset of~$I$ with its natural
$\Gal(k^\sep/k)$-action. For $J \subset I$ we obtain the usual
notion of parabolics of type~$J$. The section of~$\Dscr$
corresponding to the empty subset of~$I$ is defined over~$k$ and
$\Pscr_{\emptyset}$ is the scheme of Borel subgroups of~$G$.

We denote by $\Puni \to \Par$ the universal parabolic subgroup and
for $J \in \Dscr(S)$ we write $\Puni_J$ for its pullback to~$\Par_J \air
\Par$. We denote by~$\Uuni_J$ the unipotent radical
of~$\Puni_J$.

\secstart{}\label{\relposition} Let $S$ be a $k$-scheme. Let $P$
and~$Q$ two parabolic subgroups of~$G_S$. We say that $P$ and $Q$
are in {\it standard position\/} if the following conditions hold
(which are mutually equivalent by [\SGAthree], XXVI,~4.5):
\assertionlist
\assertionitem The intersection $P \cap Q$ is smooth.
\assertionitem Locally for the Zariski topology on the basis, $P \cap Q$
contains a maximal torus of~$G$.
\assertionitem Locally for the fpqc-topology on the basis, $P \cap Q$
contains a maximal torus of~$G$.

Let $\SP$ be the subfunctor of $\Par \times \Par$ of pairs $(P,Q)$
that are in standard position. By loc.\ cit., $\SP$ is
representable by a smooth quasi-projective scheme over~$k$. If $S$
is the spectrum of a field, any two parabolics of $G_S$ are in standard
position. Hence the monomorphism $\SP \to \Par \times \Par$ is
bijective. In fact, it can be shown that $\SP$ is the disjoint
union of the $G$-orbits in $\Par \times \Par$, in the
scheme-theoretic sense. For sections $J$, $K \in \Dscr(S)$ we
denote by $\SP_{J,K}$ the inverse image of $\Par_J \times \Par_K$
in~$\SP_S$.

The group $G$ acts on~$\SP$ by simultaneous conjugation. The fppf
quotient $G\backslash \SP$ is representable by a finite \'etale
$k$-scheme~$\RP$. Let
$$
\rbf\colon\SP \arr \RP
$$
be the canonical morphism. There exists a unique surjective
morphism of finite \'etale $k$-schemes $q\colon \RP \ar \Dscr \times
\Dscr$ such that the diagram
$$
\matrix{\SP & \arvarover(25){\rbf} & \RP \cr
\add && \addright{q} \cr
\Par \times \Par & \arvarover(25){\tbf \times \tbf} & \Dscr \times
\Dscr\cr}
$$
is commutative.

On $k^\sep$-valued points we have
$$
\RP(k^\sep) \cong \coprod_{J,K \subset I} {}^JW^K
$$
as sets with $\Gal(k^\sep/k)$-action, and ${}^JW^K = q^{-1}(J,K)$.
Hence we obtain a morphism $\iota\colon \RP \to W$ whose
restriction to $q^{-1}(J,K)$ is the inclusion ${}^JW^K \air W$. We
set
$$
\relpos := \iota \circ \rbf \colon \SP \arr W\, .
$$
Whenever we write $\relpos(P,Q)$ it shall be understood that $P$
and~$Q$ are in standard position. For $x \in W$ we define $\SP^x
:= \relpos^{-1}(x)$.

\secstart{}\label{\relposoverkbar} Over $k^\sep$ (or any other
separably closed extension of~$k$) we can describe the morphism
``$\relpos$'' as follows: We use the canonical isomorphism of the
Weyl group of the abstract root datum of~$G$ with the set of
$G(k^\sep)$-orbits in $\Par_{\emptyset}(k^\sep) \times
\Par_{\emptyset}(k^\sep)$. For $(B,B') \in
\Par_{\emptyset}(k^\sep) \times \Par_{\emptyset}(k^\sep)$ we
denote by $\relpos(B,B') \in W$ the corresponding
$G(k^\sep)$-orbit.

Now let $J$ and~$K$ be arbitrary subsets of~$I$. For $P \in
\Par_J(k^\sep)$ and $Q \in \Par_K(k^\sep)$ the relative
position $\relpos(P,Q) \in {}^JW^K$ is the unique
minimal element (with respect to the Bruhat order) in the set
$$
\bigl\{\relpos(B,B') \bigm| B \subset P, B' \subset Q\bigr\}\, .
$$
The map $(P,Q) \asr \relpos(P,Q)$ gives a bijection between the
set of $G(k^\sep)$-orbits in $\Par_J(k^\sep) \times
\Par_K(k^\sep)$ and the set~${}^JW^K$.

Alternatively we can compute $\relpos(P,Q)$ as follows: Choose a
maximal torus~$T$ which is contained in $P \cap Q$. The choice
of~$T$ provides an identification of~$W$ with $N_G(T)/T$. There
exists an $n \in N_G(T)$ such that $P$ and~$n(Q)$ contain a common
Borel subgroup and the class of~$n$ in $W_J\backslash W/W_K$
depends only on $(P,Q)$. Its unique representative in ${}^JW^K$ is
equal to $\relpos(P,Q)$.

\secstart{}\label{\refineparabolics} For $(P,Q) \in \SP(S)$ define
{\it the refinement of~$P$ with respect to~$Q$\/} to be
$$
\Ref_Q(P) := (P \cap Q)U_P = U_P(P \cap Q)\, .
$$
This is again a parabolic subgroup of~$G$ whose unipotent radical
is $U_P(P \cap U_Q) = (P \cap U_Q)U_P$. Indeed, it suffices to
show this locally for the fpqc topology hence we can assume that
$P \cap Q$ contains a split maximal torus. Then the proof is the
same as in [\BorTits],~4.4.

We refer to (2.13) for the description of $\Ref_Q(P)$ in the case that
$P$ and $Q$ are parabolics of $\GL_N$.

Suppose $P \in \Pscr_J(S)$ and $Q \in \Pscr_K(S)$ are in standard
position, with $\relpos(P,Q) = w \in {}^JW^K$. Then $\Ref_Q(P)$ is
of type $J \cap {}^wK$.

\secstart{}\label{\goodposition} Let $P$ and $Q$ be two parabolic
subgroups of~$G_S$. We say that $P$ and $Q$ are {\it in good
position\/} if the following equivalent assertions hold:
\assertionlist
\assertionitem Zariski-locally on the basis, $P$
and $Q$ contain a common Levi subgroup.
\assertionitem
fpqc-locally on the basis, $P$ and $Q$ contain a common Levi
subgroup.
\assertionitem $P$ and $Q$ are in standard position and
for every geometric point $\sbar$ of~$S$ we have that $P_{\sbar}$
and $Q_{\sbar}$ contain a common Levi subgroup.
\assertionitem $P$
and $Q$ are in standard position and for every geometric point
$\sbar$ of~$S$ we have $J_{\sbar} = {}^{w_{\sbar}}(K_{\sbar})$,
where $J_{\sbar}$ and $K_{\sbar}$ are the types of $P_{\sbar}$ and
$Q_{\sbar}$, respectively, and where $w_{\sbar} =
\relpos(P_{\sbar},Q_{\sbar})$.

(This corresponds to what in [\BMGSAS], section~3, was called ``in
optimal position''.)

\secstart{Lemma}:\label{\Levigood} {\sl Let $J$, $K \subset I$ be
sets of simple roots and let $x \in {}^KW^J$ be such that $K =
{}^xJ$. Let $Q$ be a parabolic subgroup of~$G_S$ of type~$K$ and
let $M$ be a Levi subgroup of~$Q$. Then there exists a unique
parabolic subgroup $P$ of~$G_S$ of type~$J$ such that $M$ is a
common Levi subgroup of $P$ and~$Q$ and such that $\relpos(Q,P) =
x$. (In particular, $P$ and~$Q$ are then in good position).}

\proof: Let $P$ and $P'$ be two parabolics of type~$J$ such that
$\relpos(Q,P) = \relpos(Q,P') = x$ and such that $M \subset P \cap
P'$. Then it follows from [\LusztPar],~8.4, that $\relpos(P,P') = 1$ and
hence $P = P'$. This proves the unicity.

We omit the proof of the existence as we will not need this in the
sequel.

\secstart{Lemma}:\label{\Intersectgood} {\sl Let $P$ and~$Q$ be
two parabolics of~$G_S$ which are in good position. Then we have
$U_P \cap Q = U_P \cap U_Q$.}

\proof: The question is local on~$S$ for the fpqc topology; hence
we can assume that there exists a common Levi subgroup of $P$
and~$Q$ whose connected center is a split torus. Now the proof is
the same as in [\LusztPar],~8.6.

\secstart{} Let $P \in \Par_J(S)$ and $Q \in \Par_K(S)$ with $K =
{}^{w_0}J$  where $w_0$ is the longest element of~$W$. Let $x$ be
the element of minimal length in the double coset $W_Jw_0W_K$.
Then $P$ and~$Q$ are in opposition (i.e., $P \cap Q$ is a
common Levi subgroup of $P$ and~$Q$), if and only if $\relpos(P,Q)
= x$.

\secstart{}\label{\parcorrespondence} Let $P$ and~$Q$ be two
parabolics of $G_S$ which are in good position. Then every
parabolic subgroup $P'$ of~$P$ is in standard position with~$Q$
and we have $\relpos(P',Q) = \relpos(P,Q)$. The maps $P' \asr
\Ref_{P'}(Q)$ and $Q' \asr \Ref_{Q'}(P)$ define mutually inverse
bijections
$$
\bigl\{\hbox{parabolic subgroups of $P$}\bigr\} \longleftrightarrow
\bigl\{\hbox{parabolic subgroups of $Q$}\bigr\}\, .
$$
Moreover, $P'$ and $\Ref_{P'}(Q)$ are in good position and we have
$$
\relpos\bigl(P',\Ref_{P'}(Q)\bigr) = \relpos(P,Q)\, .
$$
In particular we see that $\Ref_P(Q) = Q$.

\secstart{Example}:\label{\GLExa} Let $G = \GL_n$. Then $G$ is
split over~$k$. Associating to a flag in~$\Oscr^n_S$ its
stabilizer defines an isomorphism between the scheme of flags and
the scheme~$\Par$. We use this isomorphism to identify flags in
$\Oscr^n_S$ and parabolics of~$G_S$.

The Weyl group $W$ can be identified with $S_n$ such that $I$ is
the set of transpositions $\tau_{\alpha} = (\alpha\; \alpha+1)$
for $\alpha = 1,\dots, n-1$. If $\Gamma = (\Gamma^i)$ is a flag
such that all $\Gamma^i$ have constant rank, its type $J \subset
I$ is determined by the rule that $\tau_\alpha \notin J$ if and
only if there exists an index~$i$ with $\rk_{\Oscr_S}(\Gamma^i) =
\alpha$.

Let $\Gamma = (\Gamma^i)$ and $\Delta = (\Delta^j)$ be two flags
in $\Oscr_S^n$. Then the following conditions are equivalent:
\assertionlist \assertionitem The parabolics associated to
$\Gamma$ and $\Delta$ are in standard position. \assertionitem For
all $i$ and $j$, the submodule $\Gamma^i + \Delta^j \subset
\Oscr_S^n$ is locally a direct summand. \assertionitem
Zariski-locally on $S$ there exists a basis $\{e_1,\dots,e_n\}$
of~$\Oscr_S^n$, such that for all $i$ and $j$ there exists a
subset $I_{i,j}$ of $\{1,\dots,n\}$ with $\Gamma^i + \Delta^j =
\bigoplus_{\alpha \in I_{i,j}}\Oscr_S\cdot e_{\alpha}$.

If these conditions are satisfied, the relative position of
$\Gamma$ and~$\Delta$ is completely determined by the function
$(i,j) \mapsto \rk_{\Oscr_S}(\Gamma^i + \Delta^j)$.

As an example, for $J$, $K \subset I$, let $x$ be the element of
minimal length in $W_Jw_0W_K$, where $w_0$ is the longest element
in~$W$. Let $\Gamma$ and~$\Delta$ be flags of types $J$ and~$K$,
respectively, which are in standard position. Then we have
$$
\relpos(\Gamma,\Delta) = 1
\quad\Longleftrightarrow\quad
\rk_{\Oscr_S}(\Gamma^i + \Delta^j) =
\max\bigl(\rk_{\Oscr_S}(\Gamma^i),\rk_{\Oscr_S}(\Delta^j)\bigr)
\quad \hbox{for all $i$, $j$}\, ,
$$
and
$$
\relpos(\Gamma,\Delta) = x
\quad\Longleftrightarrow\quad
\rk_{\Oscr_S}(\Gamma^i + \Delta^j) = \min\bigl(n,\rk_{\Oscr_S}(\Gamma^i) +
\rk_{\Oscr_S}(\Delta^j)\bigr)
\quad \hbox{for all $i$, $j$}\, .
$$

If $\Gamma$ and $\Delta$ are flags in standard position with
stabilizers $P$ and~$Q$, respectively, the flag corresponding to
$\Ref_Q(P)$ is given by the collection of submodules
$(\Gamma^{i-1} \cap \Delta^j) + \Gamma^i$ for all $i$ and~$j$. This is a
refinement of the flag $\Gamma$.

If $\Gamma$ is a flag with associated parabolic $P$, the choice of
a Levi subgroup of~$P$ corresponds to the choice of a
decomposition $\Oscr_S^n = \bigoplus_{j=1}^r M_j$ such that
$\Gamma^i = \bigoplus_{j > i} M_{\pi(j)}$ for some permutation
$\pi \in S_r$.


\paragraph{A semi-linear variation on a theme of Lusztig}

\nnsection{}In this section we consider a reductive group $G$
over~$\Fq$. As in Lusztig's paper~[\LusztPar], we define, for $J$ a set
of simple reflections in the Weyl group, a variety~$Z_J$ equipped
with an action of~$G$. This is a semi-linear variant of the
variety examined by Lusztig. The main result of this section,
Theorem~(3.25), concerns a classification of the $G$-orbits
in~$Z_J$. This result will be used in the next section to prove
our main classification theorem for $F$-zips.

Throughout this section, $q$ is a fixed power of a prime
number~$p$, and $S$ is a scheme over~$\F_q$. Note that in~[\LusztPar],
Lusztig writes $P^Q$ for what we call~$\Ref_Q(P)$.

\secstart{}\label{\GhatSetup} Let $\Ghat$ be a possibly
disconnected reductive group over~$\Fq$ and denote by~$G$ its
identity component. We keep the notations of~\groupnotation; note
that $G$ is indeed quasi-split. Further we fix a connected
component $G^1$ of~$\Ghat$. Let $\Fdbar$ be an algebraic closure
of~$\F_q$ and denote by $\sigma\colon x \asr x^q$ the arithmetic
Frobenius in $\Gal(\Fdbar/\Fq)$. It acts on~$(W,I)$.

There is a unique $\Gal(\Fdbar/\Fq)$-equivariant isomorphism
$\delta\colon (W,I) \to (W,I)$ of Coxeter systems such that for
all $g \in G^1(\Fdbar)$ and $P \in \Pscr_J(\Fdbar)$ we have ${}^gP
\in \Pscr_{\delta(J)}(\Fdbar)$.

If there is no risk of confusion we simply write $F\colon \Ghat
\to \Ghat$ for the morphism $\Fr_{\Ghat}\colon \Ghat \to
\Ghat^{(q)} = \Ghat$ that was defined in~\Frdef. It is an
endomorphism of~$\Ghat$.

\secstart{}\label{\JKx} We fix the following data: Let $J$ and $K$
be subsets of~$I$ and $x \in W$ such that ${}^x\delta(J) = K$ and
$x \in {}^KW^{\delta(J)}$. We assume that $J$ and~$x$ (hence
also~$K$) are defined over~$\Fq$, i.e., $\sigma(J) = J$ and
$\sigma(x) = x$.

\secstart{}\label{\DefineZtilde} Let $\tilde{Z}_J$ be the
$\Fq$-scheme given by the cartesian square
$$
\matrix{
\tilde{Z}_J & \arr & \SP^x\cr
\add & & \add\cr
\Par_J \times \Par_K \times G^1 & \arrover{f} & \Par_K \times
\Par_{\delta(J)}\cr
}
$$
where the morphism $f$ is given on points by $(P,Q,g) \mapsto
\bigl(Q,{}^gF(P)\bigr)$. If $S$ is an $\Fq$-scheme then the
$S$-valued points of $\tilde{Z}_J$ are the triples $(P,Q,g)$ with
$P$ and~$Q$ parabolics of~$G_S$ of types $J$ and~$K$,
respectively, and with $g \in G^1(S)$ an element such that $Q$
and~${}^gF(P)$ are in relative position~$x$. In particular, $Q$
and~${}^gF(P)$ are then in good position; see~\goodposition. The
forgetful morphism $(P,Q,g) \asr (P,Q)$ makes $\tilde{Z}_J$ into a
scheme over $\Par_J \times \Par_K$.

We define an action of~$G$ on $\Ztilde_J$ given on $S$-valued
points by
$$
h \cdot (P,Q,g) := \bigl({}^hP, {}^hQ, hgF(h)^{-1}\bigr)\, .
$$
It is easily seen that this is well-defined.

\secstart{} For $u \in {}^JW^K$, let $\tilde{Z}^u_J$ be the
subscheme of~$\tilde{Z}_J$ of triples $(P,Q,g)$ with $\relpos(P,Q)
= u$. The natural morphism
$$
\coprod_{u \in {}^JW^K} \tilde{Z}^u_J \longrightarrow \tilde{Z}_J
$$
is a bijective monomorphism.

Fix $u \in {}^JW^K$. Let $L :=  \delta\bigl(J \cap
{}^{ux}\delta(J)\bigr) = \delta(J \cap {}^uK)$ and consider the
morphism $h\colon \tilde{Z}^u_J \to \Par_K \times \Par_L$ given on
points by $(P,Q,g) \mapsto \bigl(Q,{}^gF(\Ref_{Q}(P))\bigr)$.
Define a scheme $\Ytilde^u_J$ by the fibre product diagram
$$
\matrix{
\Ytilde^u_J & \arr & \SP_{K,L}\cr
\add & & \add\cr
\tilde{Z}^u_J & \arrover{h} & \Par_K \times \Par_L\, .\cr
}
$$
On points this means that we are considering triples $(P,Q,g)$ in
$\tilde{Z}^u_J$ with the additional requirement that $Q$ and
${}^gF\bigl(\Ref_{Q}(P)\bigr)$ are in standard position.

Note that the $G$-action on~$\Ztilde_J$ preserves $\Ztilde_J^u$
and~$\Ytilde_J^u$.

\secstart{}
Let
$$
J_1 := J \cap {}^{ux}\delta(J) = J \cap {}^uK
\quad\hbox{and}\quad
K_1 := {}^x \delta(J_1)\, .
$$
Define a morphism $\tilde\vartheta \colon \Ytilde_J^u \to
\tilde{Z}_{J_1}$ by $\tilde\vartheta(P,Q,g) = (P_1,Q_1,g)$ with
$$
P_1 := \Ref_{Q}(P)
\quad\hbox{and}\quad
Q_1 := \Ref_{{}^gF(\Ref_{Q}(P))}(Q) = \Ref_{{}^gF(P_1)}(Q)\, .
$$
To see that $\tilde\vartheta$ is well-defined, we need to check
that $P_1$ and $Q_1$ are parabolics of types $J_1$ and~$K_1$,
respectively, and that $\relpos\bigl(Q_1,{}^gF(P_1)\bigr) = x$.
That $P_1$ has type~$J_1$ is immediate from~\refineparabolics.
Next remark that ${}^gF(P_1) \subset {}^gF(P)$, so
by~\parcorrespondence\ we have $\relpos\bigl(Q,{}^gF(P_1)\bigr) =
x$. Again using~\refineparabolics\ we then easily verify that
$Q_1$ has type~$K_1$, and by~\parcorrespondence\ we conclude that
$\relpos\bigl(Q_1,{}^g F(P_1)\bigr) = x$.

\secstart{}\label{\usequences} Consider a sequence $\ubf =
(u_0,u_1,\ldots)$ of elements of~$W$. Define a sequence of subsets
$J_n \subset I$ by setting $J_0 := J$ and $J_{n+1} := J_n \cap
{}^{u_nx}\delta(J_n)$. Set $K_n := {}^x \delta(J_n)$. Let
$\Tscr(J)$ be the set of sequences~$\ubf = (u_0,u_1,\ldots)$ such that
for all $n \geq 0$ we have
$$
u_n \in {}^{J_n}W^{K_n}
\quad\hbox{and}\quad
u_{n+1} \in W_{J_{n+1}}u_n W_{K_n}\, .\lformno\sublabel{\Combinatcond}
$$
(These conditions imply that in fact $u_{n+1} \in u_n W_{K_n}$.)  
By construction, $J_{n+1} \subseteq J_n$ and $K_{n+1} \subseteq K_n$ for
all~$n$. Hence there exists an index~$N$ such that $J_{n+1} = J_n$ and
$K_{n+1} = K_n$ for all $n \geq N$. Writing $J_\infty := J_n$ and $K_\infty
:= K_n$ for $n \gg 0$, we find that $J_\infty = {}^{u_n} K_\infty$ for $n
\geq N$. If $\ubf \in \Tscr(J)$ and $n \geq N$ then the two conditions
in~\Combinatcond\ readily imply that $u_{n+1} = u_n$. Set $u_\infty := u_n$
for any $n \geq N$.

\secstart{Lemma}:\label{\Combinatlemma} {\sl The map $\Tscr(J) \to
W$ defined by $\ubf \mapsto u_\infty$ gives a bijection $\Tscr(J)
\alr {}^JW$.}

\proof: Set $\Jtilde = K_0$ and $\Jtilde' = J$ and let
$\varepsilon$ be the automorphism $w \asr \delta(x^{-1}wx)$ of
$W$. Then the set $\Tscr(J)$ is nothing but Lusztig's set
$\Tscr(\Jtilde,\varepsilon)$ as defined in [\LusztPar],~2.2, and our
claim follows from [\LusztPar],~2.5.

\secstart{} Let $\ubf = (u_0,u_1,\ldots) \in \Tscr(J)$. Let
$N(\ubf)$ be the smallest non-negative integer such that $J_{n+1}
= J_n$ for all $n \geq N(\ubf)$; as we have seen this implies that
also $K_{n+1} = K_n$ and $u_{n+1} = u_n$ for all $n \geq N(\ubf)$.

For $r \geq 0$ we write $\ubf_r := (u_r,u_{r+1},\ldots)$, which is an
element
of~$\Tscr(J_r)$. In particular, $\ubf = \ubf_0$. Note that $N(\ubf_r) =
\max\{0,N(\ubf_0)-r\}$.

By induction on $N(\ubf)$ we now define schemes $\Ytilde_J^\ubf$
together with morphisms $\Ytilde_J^\ubf \to \Ytilde_J^{u_0}$. If
$N(\ubf) = 0$ then we set $\Ytilde_J^\ubf := \Ytilde_J^{u_0}$,
mapping identically to itself. Next  assume that $N(\ubf) = N$ and
that for all $L \subset I$ and $\vbf \in \Tscr(L)$ with $N(\vbf) <
N$ the morphism of schemes $\Ytilde_L^\vbf \to \Ytilde_L^{v_0}$
has been defined. Then we define $\Ytilde_J^\ubf$ by the fibre
product diagram
$$
\matrix{
\Ytilde_J^\ubf & \arvar(25) & \Ytilde^{u_0}_J\cr
\add & &  \addright{\tilde\vartheta} \cr
\Ytilde_{J_1}^{\ubf_1} & \ar \Ytilde_{J_1}^{u_1} \ar &
\tilde{Z}_{J_1}\rlap{\ .}\cr
}
$$

On points this means the following. If $N(\ubf) = 0$ then $\ubf=
(u,u,\ldots)$ is a constant sequence, and we just consider the
scheme $\Ytilde_J^u$. Next suppose $N(\ubf) = 1$, which means that
$\ubf = (u_0,u_1,u_1,\ldots)$ for some $u_0 \neq u_1$. In this
case, the points of $\Ytilde_J^\ubf$ are the points $(P,Q,g)$ of
$\Ytilde_J^{u_0}$ such that the associated triple $(P_1,Q_1,g) :=
\tilde\vartheta(P,Q,g)$ lies in $\Ytilde_{J_1}^{u_1}
\hookrightarrow \tilde{Z}_{J_1}$. In general we have a diagram
$$\matrix{
& & & & \Ytilde_J^{u_0} & \airr & \tilde{Z}_J \cr & & & &
\addright{\tilde\vartheta} \cr & & \Ytilde_{J_1}^{u_1} & \airr &
\tilde{Z}_{J_1} \cr & & \addright{\tilde\vartheta} \cr
\Ytilde_{J_2}^{u_2} & \airr & \tilde{Z}_{J_2} \cr
\addright{\tilde\vartheta} \cr \cdots \cr }$$ and the points of
$\Ytilde_J^\ubf$ are those triples $(P,Q,g)$ in $\Ytilde_J^{u_0}$
that under each subsequent map $\tilde\vartheta$ land inside
$\Ytilde_{J_n}^{u_n} \hookrightarrow \tilde{Z}_{J_n}$.

Note that the map $\Ytilde^{\ubf}_J \to \Ytilde^{u_0}_J$ is a
monomorphism and that the $G$-action on $\Ytilde^{u_0}_J$
preserves~$\Ytilde^{\ubf}_J$.

\secstart{}\label{\defineZ}
Let $\ubf = (u_0,u_1,\dots) \in \Tscr(J)$ and set $u = u_0$.  The schemes
$$
\Ytilde_J^\ubf \airr \Ytilde_J^u \airr \tilde{Z}_J^u
\airr \tilde{Z}_J\lformno\sublabel{\tildeinclusion}
$$
are schemes over $\Par_J \times \Par_K$. Recall that we denote by
$\Puni_J$ the universal parabolic group scheme over~$\Par_J$ and
by $\Uuni_J$ its unipotent radical. Then $F(\Puni_J)$ is again a
parabolic subgroup scheme of $G \times \Par_J$ over~$\Par_J$,
which has $F(\Uuni_J)$ as its unipotent radical. (In fact, as $J$
is defined over~$\Fq$, so is $\Par_J$, and $F(\Uuni_J)$ is none
other than the pull-back of~$\Uuni_J$ via the morphism
$F_{\Par_J}\colon \Par_J \to \Par_J$.)

Write $\Uuni_{J,K}$ for the scheme $F(\Uuni_J) \times \Uuni_K$,
but with a new group scheme structure given on points by
$(v_1,v_2) \cdot (v_1^\prime,v_2^\prime) = (v_1^\prime
v_1,v_2v_2^\prime)$. Then $\Uuni_{J,K}$ acts from the left on all
four schemes in~\tildeinclusion\ by
$$
(v_1,v_2) \cdot (P,Q,g) = (P,Q,v_2g v_1)\, .
$$

We define
$$
Y_J^\ubf \airr Y_J^u \airr Z_J^u \airr Z_J\lformno\sublabel{\inclusion}
$$
to be the fppf quotient sheaves of the schemes in~\tildeinclusion\
by this action of~$\Uuni_{J,K}$. More informally we could write
$Z_J = \Uuni_K\backslash\tilde{Z}_J/F(\Uuni_J)$, and similarly for
the other quotients. If $(P,Q,g) \in \tilde{Z}_J(S)$ then we write
$[P,Q,g]$ for its image in~$Z_J(S)$.

It readily follows from the definitions that the $G$-action on
$\Ztilde_J$ induces a $G$-action on $Z_J$, and hence on all other
quotients in~\inclusion. Further it follows from [\SGAthree],
Exp.~XXVI,~2.2 that for an affine scheme~$S$ the canonical
morphism $\Ztilde_J(S) \to Z_J(S)$ is surjective.

\secstart{} Our next goal is to show that the sheaves
in~\inclusion\ are representable by schemes.

The quotient $\Puni_J/\Uuni_J$ is representable by a reductive
group scheme over~$\Par_J$. Let $H$ be defined by the cartesian
diagram
$$
\matrix{
H & \arr & \Puni_J/\Uuni_J \cr
\add && \add \cr
\Par_J \times \Par_K & \arrover{\alpha} & \Par_J\, ,}
$$
with $\alpha$ given by $(P,Q) \mapsto F(P)$. For an affine
scheme~$S$, the $S$-valued points of~$H$ are given by triples
$\bigl(P,Q,yU_{F(P)}(S)\bigr)$, where $P$ and~$Q$ are parabolic
subgroups of $G_S$ of types~$J$ and~$K$, respectively, and where
$y \in F(P)\bigl(S\bigr)$.

Define a right action
$$
Z_J \times_{(\Par_J \times \Par_K)} H \arr Z_J
$$
as follows: For an affine scheme~$S$, a point $z = [P,Q,g] \in
Z_J(S)$, and $h = (P,Q,yU_{F(P)}) \in H(S)$, we set
$$
z \cdot h = [P,Q,gy] \in Z_J(S)\, .
$$

\secstart{Lemma}:\label{\Zastorsor} {\sl This action makes $Z_J$
into an $H$-torsor over $\Par_J \times \Par_K$ for the Zariski topology.}

\proof: Let $P \in \Par_J(S)$ and $Q \in \Par_K(S)$ and suppose we
have $g$, $g' \in G^1(S)$ such that $\relpos\bigl(Q,{}^gF(P)\bigr)
= \relpos\bigl(Q,{}^{g'}F(P)\bigr) = x$. Locally on~$S$ we can
find $b \in Q$ such that $g' \in bgF(P)$. Let $M$ be a common Levi
subgroup of ${}^gF(P)$ and~$Q$ (which we can find Zariski-locally
on~$S$, as ${}^gF(P)$ and~$Q$ are in good position). Then we have
$b = vm$ with $v \in U_Q$ and $m \in M$. As $M \subset {}^gF(P)$,
we have $g' \in vmgF(P) = vgF(P)$. This proves that the action is
transitive.

Now assume that for $g \in G^1$ with
$\relpos\bigl(Q,{}^gF(P)\bigr) = x$ there exist elements $y$, $y'
\in F(P)$ such that $gy' \in U_QgyU_{F(P)}$. Then possibly after
multiplying~$y$ from the right by an element of~$U_{F(P)}$ we may
assume that there is a $v \in U_Q$ with $gy' = vgy$. But then $v
\in U_Q \cap {}^gF(P) = U_Q \cap {}^gU_{F(P)}$, where the last
equality holds by~\Intersectgood. Hence there exists a $v' \in
U_{F(P)}$ such that $v = gv'g^{-1}$. This gives that $y' \in
U_{F(P)} \cdot y = y \cdot U_{F(P)}$, proving that the action is
free.

\secstart{Corollary}: {\sl The fppf sheaves $Y_J^\ubf$, $Y_J^u$,
$Z_J^u$ and~$Z_J$ are representable by schemes.}

\secstart{Lemma}: {\sl The morphism $\tilde\vartheta \colon
\Ytilde_J^u \to \tilde{Z}_{J_1}$ induces a morphism $\vartheta
\colon Y_J^u \to Z_{J_1}$.}

\proof: Let $(P,Q,g)$ be an $S$-valued point of~$\Ytilde_J^u$, and
let $(P_1,Q_1,g)$ be its image under~$\tilde\vartheta$. Then $P_1
\subseteq P$ and $Q_1 \subseteq Q$, so $U_{F(P)} \subseteq
U_{F(P_1)}$ and $U_Q \subseteq U_{Q_1}$. But then it is immediate
from the definitions that the composed morphism $\Ytilde_J^u \to
\tilde{Z}_{J_1} \to Z_{J_1}$ factors modulo the action
of~$\Uuni_{J,K}$.

\secstart{} Let $\ubf = (u_0,u_1,u_2,\ldots) \in \Tscr(J)$. For $n
\geq 0$ let $\ubf_n = (u_n,u_{n+1},\ldots) \in \Tscr(J_n)$. As an
immediate consequence of the definition of the schemes
$\Ytilde_J^\ubf$ and their quotients~$Y_J^\ubf$ we obtain
$G$-equivariant morphisms
$$
\tilde\vartheta\colon \Ytilde_J^{\ubf_n} \to \Ytilde_{J_1}^{\ubf_{n+1}}\, .
$$
inducing $G$-equivariant morphisms
$$
\vartheta\colon Y_J^{\ubf_n} \to Y_{J_1}^{\ubf_{n+1}}\, .
$$

\secstart{}\label{\ZJunionYJu} The natural map $\coprod_{\ubf \in 
\Tscr(J)} \tilde{Y}_J^\ubf \to \tilde{Z}_J$ is a bijective 
monomorphism. Passing to quotients modulo~$\Uuni_{J,K}$ we readily 
find that $\coprod_{\ubf \in \Tscr(J)} Y_J^\ubf \to Z_J$ is a 
bijective monomorphism, too. (Use that $\tilde{Z}_J \to Z_J$ is 
surjective on underlying topological spaces.) In particular, if $k$ 
is an algebraically closed field then we have a bijection
$$
\coprod_{\ubf \in \Tscr(J)} Y_J^\ubf(k) \arriso Z_J(k)\, .
$$

Our main goal for the rest of this section is to show that the 
$G$-action on the schemes $Y_J^\ubf$ is transitive. Along the way we
shall also compute the dimension of the schemes~$Y_J^\ubf$.

\secstart{Lemma}:\label{\thetasurjective} {\sl The morphism
$\tilde{\vartheta}\colon \Ytilde_{J_n}^{\ubf_n} \arr
\Ytilde_{J_{n+1}}^{\ubf_{n+1}}$ is an isomorphism.}

\proof: Without loss of generality we can assume $n = 0$. Suppose
$\tilde{\vartheta}(P,Q,g) = \tilde{\vartheta}(P',Q',g') =:
(P_1,Q_1,g_1)$. Clearly, $g = g_1 = g'$. Further, $P$ and~$P'$ are
parabolics of the same type, and they both contain~$P_1$; hence $P
= P'$. The same argument shows that $Q = Q'$. Hence
$\tilde{\vartheta}$ is a monomorphism.

Now let $(P_1,Q_1,g)$ be an $S$-valued point
of~$\Ytilde_{J_1}^{\ubf_1}$. Let $P$ be the unique parabolic of
type~$J$ that contains~$P_1$, and let $Q$ be the unique parabolic
of type~$K$ containing~$Q_1$. These exist by [\SGAthree],
Exp.~XXVI,~3.8. Then $(P_1,Q_1)$ is in standard position, $P
\supset P_1$, and $Q \supset Q_1$; hence $(P,Q)$ is in standard
position, too. In a similar way we see that
$\bigl({}^gF(\Ref_{Q}(P)),Q\bigr)$ is in standard position.

By definition of~$\Tscr(J)$ we have $\relpos(P_1,Q_1) \in uW_K$
with $u = u_0 \in {}^JW^K$. It follows that $\relpos(P,Q) = u$.
Similarly, as $x\in {}^KW^{\delta(J)}$ and
$\relpos\bigl(Q_1,{}^gF(P_1)\bigr) = x$, we also have
$\relpos\bigl(Q,{}^gF(P)\bigr) = x$. Hence it remains to see that
$\Ref_{Q}(P) = P_1$ and $\Ref_{{}^gF(P_1)}(Q) = Q_1$. For this we
may work fppf-locally on~$S$.

We write $\relpos\bigl(P_1,Q_1) = uw$ with $w \in W_K$. Working
fppf-locally we can assume that $G_S$ is split and hence we can
find Borel subgroups $B \subset P_1$ and $C \subset Q_1$ such that
$\relpos(B,C) = uw$. As we have $\ell(uw) = \ell(u) + \ell(w)$,
there exists a Borel subgroup $D$ of~$G_S$ such that $\relpos(B,D)
= u$ and $\relpos(D,C) = w$. As $C \subset Q_1 \subset Q$ and $w
\in W_K$, we see that $D \subset Q$. As $B \subset P$ and $D
\subset Q$ and $\relpos(B,D) = \relpos(P,Q)$, we have $B \subset
\Ref_{Q}(P)$. But then $P_1$ and $\Ref_{Q}(P)$ have the same type
and have a Borel subgroup in common; hence they are equal.

It remains to be shown that $Q_1 = \Ref_{{}^gF(P_1)}(Q)$. Set
$Q_1' := \Ref_{{}^gF(P_1)}(Q)$. By \parcorrespondence\ we have
$\relpos\bigl(Q_1,{}^gF(P_1)\bigr) = x =
\relpos\bigl(Q_1',{}^gF(P_1)\bigr)$, so there exists an $h \in
{}^gF(P_1)$ with $Q_1 = {}^hQ'_1$. Moreover we have ${}^hQ \supset
{}^hQ'_1 = Q_1$ and $Q \supset Q_1$ and hence ${}^hQ = Q$.
Therefore,
$$
Q_1 = {}^hQ'_1 = {}^h\Ref_{{}^gF(P_1)}(Q) = \Ref_{{}^gF(P_1)}(Q)\, ,
$$
and the proof is complete.

\secstart{Lemma}:\label{\thetabijective} {\sl The morphism
$\vartheta\colon Y_{J_n}^{\ubf_n} \arr Y_{J_{n+1}}^{\ubf_{n+1}}$
induces an isomorphism of fppf quotient sheaves
$$
\bar{\vartheta}\colon G\backslash Y_{J_n}^{\ubf_n} \arriso G\backslash
Y_{J_{n+1}}^{\ubf_{n+1}}.
$$}

\proof: Without loss of generality we can assume $n = 0$. We
consider $S$-valued points, where $S$ is an affine $\Fq$-scheme.
Note that the quotient map $\Ytilde^{\ubf}_J(S) \to Y^\ubf_J(S)$
is surjective by [\SGAthree], Exp.~XXVI,~2.2. Further let us recall
that for $(P,Q,g) \in \Ytilde^{\ubf}_J(S)$ we denote by $[P,Q,g]$
its image in~$Y^\ubf_J$.

By \thetasurjective\ we only have to show that $\bar{\vartheta}$
is a monomorphism. Let $[P,Q,g]$ and $[P',Q',g']$ be two
$S$-valued points of~$Y^\ubf_J$ such that
$\vartheta\bigl([P,Q,g]\bigr)$ and
$\vartheta\bigl([P',Q',g']\bigr)$ are in the same $G(S)$-orbit. We
want to show that $[P,Q,g]$ and $[P',Q',g']$ are fppf-locally in
the same $G(S)$-orbit. As $\vartheta$ is $G$-equivariant, we can
assume that $\vartheta([P,Q,g]) = \vartheta([P',Q',g']) =:
[P_1,Q_1,g_1]$.

As $P$ and~$P'$ have the same type and both contain~$P_1$, we get
$P = P'$. Similarly, $Q = Q'$.

By definition, $\relpos\bigl(Q,{}^gF(P)\bigr) =
\relpos\bigl(Q,{}^{g^\prime}F(P)\bigr) = x$; hence ${}^gF(P)$ and
${}^{g^\prime}F(P)$ are both in good position to~$Q$. Let $L$
(resp.~$L'$) be a common Levi subgroup of ${}^gF(P)$ and~$Q$
(resp.~of ${}^{g^\prime}F(P)$ and~$Q$). There exists a unique $\xi
\in U_Q$ with ${}^\xi L = L'$ ([\SGAthree], Exp.~XXVI, 1.8). By
\Levigood\ this implies that ${}^{\xi g}F(P) = {}^{g^\prime}F(P)$,
and therefore $\xi g \in g'F(P)$. We can replace $g$ by~$\xi g$
and write $g' = gy$ with $y \in F(P)$. In particular, we now have
${}^gF(P) = {}^{g^\prime}F(P)$.

We have
$$
\Ref_{{}^gF(\Ref_Q(P))}(Q) = Q_1 = \Ref_{{}^{g\prime}F(\Ref_Q(P))}(Q)
$$
and this is a parabolic subgroup of~$Q$. As $Q$ and ${}^gF(P) =
{}^{g^\prime}F(P)$ are in good position, \parcorrespondence\
implies that ${}^gF\bigl(\Ref_Q(P)\bigr) =
{}^{g\prime}F\bigl(\Ref_Q(P)\bigr)$; in other words
$$
{}^gF(P_1) = {}^{g\prime}F(P_1)\, .
$$

By hypothesis, $g' \in U_{Q_1}g U_{F(P_1)}$, at least
fppf-locally. Hence can write $g' = vgu$ with $v \in U_{Q_1} =
(U_{{}^gF(P_1)} \cap Q)U_Q$ and $u \in U_{F(P_1)}$. Changing $g$
on the left by an element of~$U_Q$, we may assume that $v \in
U_{{}^gF(P_1)} \cap Q$. Write $v = gu' g^{-1}$ with $u' \in
U_{F(P_1)}$ and replace $u$ by~$u'u$. Then we have $g' = gu$ with
$u \in U_{F(P_1)} = U_{F(P)} \bigl(F(P) \cap U_{F(Q)} \bigr)$;
see~\refineparabolics.

Write $u = u_1u_2$ with $u_1 \in U_{F(P)}$ and $u_2 \in F(P) \cap
U_{F(Q)}$. Replacing $g$ by~$gu_1$ and $u$ by~$u_2$ we can further
assume that
$$
g' = gu\, , \quad \hbox{with\ } u \in F(P) \cap U_{F(Q)}\, .
\lformno\sublabel{\kappasurjective}
$$
Note that we did not use the $G$-action so far.

To finish the proof, all we now have to remark is that
fppf-locally on~$S$ we can write $u = F(v)$ with $v \in P \cap
U_Q$ (as $F\colon P \cap U_Q \to F(P) \cap U_{F(Q)}$ is an
epimorphism of fppf sheaves), and then
$$
v^{-1}\cdot [P,Q,g] =
\bigl[{}^{v^{-1}}P,{}^{v^{-1}}Q,v^{-1}gF(v)\bigr] = [P,Q,v^{-1}g']
= [P,Q,g']\, .
$$
Hence $\bar{\vartheta}$ is indeed injective.

\secstart{Lemma}:\label{\describefibre} {\sl Let $S$ be an affine
scheme. For $[P_1,Q_1,g_1] \in Y^{\ubf_1}_{J_1}(S)$, choose
$[P,Q,g] \in Y^\ubf_J(S)$ with $\vartheta\bigl([P,Q,g]\bigr) =
[P_1,Q_1,g_1]$. (This is possible by \thetasurjective\ and the
surjectivity of the map $\Ytilde^{\ubf_1}_{J_1}(S) \to
Y^{\ubf_1}_{J_1}(S)$.) Then we have a well-defined morphism
$$
\kappa\colon F(P) \cap U_{F(Q)} \arr
\vartheta^{-1}\bigl([P_1,Q_1,g_1]\bigr)
$$
given on points by $v \asr [P,Q,gv]$, and this induces an
isomorphism
$$
\bigl(F(P) \cap U_{F(Q)}\bigr) \bigm/ \bigl(U_{F(P)} \cap U_{F(Q)}
\bigr) \arriso \vartheta^{-1}\bigl([P_1,Q_1,g_1]\bigr)\, .
$$}

\proof: It is easy to check that $\kappa$ is well-defined. The
arguments of~\thetabijective, resulting in the
relation~\kappasurjective, show that $\kappa$ is an epimorphism of
sheaves.

It is clear that if $v = v'y$ for some $y \in U_{F(P)}$ then
$\kappa(v) = \kappa(v')$. Conversely, assume that $\kappa(v) =
\kappa(v')$. Then we have
$$
gv' \in U_QgvU_{F(P)} = U_Q g U_{F(P)} v\, ,
$$
so we may write $gv' = w g u v$ with $w \in U_Q$ and $u \in
U_{F(P)}$. But then $w = gv'v^{-1}u^{-1}g^{-1} \in U_Q \cap
{}^gF(P)$, so
$$
wgug^{-1} \in \bigl(U_Q \cap {}^gF(P)\bigr) {}^gU_{F(P)} =
U_{\Ref_Q({}^gF(P))} = U_{{}^gF(P)}\, ,
$$
where the last equality holds because $Q$ and~${}^gF(P)$ are in
good position. It follows that $wgu = gy$ for some $y \in
U_{F(P)}$; hence $v' = yv \in U_{F(P)} \cdot v = v \cdot
U_{F(P)}$.

\secstart{}\label{\Seedimension} Let $\Puni_J$ and $\Puni_K$ be
the universal parabolic subgroups over $\Par_J$ and~$\Par_K$,
respectively. Define an action of $F(\Puni_J) \cap U_{F(\Puni_K)}$
on~$\Ztilde_J$ over $\Par_J \times \Par_K$  by
$$
v\cdot (P,Q,g) = (P,Q,gv)\, .
$$

For $u \in {}^JW^K$ this action preserves $\Ytilde^u_J
\hookrightarrow \tilde{Z}_J$. Moreover,
$\tilde\vartheta\bigl(v\cdot(P,Q,g)\bigr) =
\tilde\vartheta(P,Q,g)$, so $F(\Puni_J) \cap U_{F(\Puni_K)}$ acts
on the fibres of~$\tilde\vartheta$. Obviously, this action
descends to an action on $Z_J$ and~$Y^u_J$. Hence for a scheme~$S$
over $\Par_J \times \Par_K$ and a section $y \in Y_J^u(S)$, we
have that $\bigl(F(\Puni_J) \cap U_{F(\Puni_K)}\bigr)_S$ acts on
the fibre $\vartheta^{-1}\bigl(\vartheta(y)\bigr)$. Now
\describefibre\ shows that
$\vartheta^{-1}\bigl(\vartheta(y)\bigr)$ is a torsor under the
affine group scheme $\bigl(F(\Puni_J) \cap U_{F(\Puni_K)}) \bigm/
\bigl(U_{F(\Puni_J)} \cap U_{F(\Puni_K)}\bigr)_S$. Moreover,
$$
\eqalign{ \bigl(F(\Puni_J) \cap U_{F(\Puni_K)}\bigr) \bigm/
\bigl(U_{F(\Puni_J)} \cap U_{F(\Puni_K)}\bigr) \arriso &
\;\bigl(F(\Puni_J) \cap U_{F(\Puni_K)}\bigr) U_{F(\Puni_J)} \bigm/
U_{F(\Puni_J)}\cr \hfill =\; & \;
U_{\Ref_{F(\Puni_K)}\bigl(F(\Puni_J)\bigr)}/U_{F(\Puni_J)}\, . }$$
Hence the dimension of the fibres of~$\vartheta$ equals
$$
\dim\bigl(F(\Uuni_{J_1})\bigr) - \dim\bigl(F(\Uuni_J)\bigr)   =
\dim(\Uuni_{J_1}) - \dim(\Uuni_J) = \dim(\Par_{J_1}) -
\dim(\Par_J)\, .
$$

Repeating this argument, we obtain a chain of morphisms
$$
Y^\ubf_J \arrover{\vartheta_1} Y^{\ubf_1}_{J_1} \arrover{\vartheta_2}
\cdots \arrover{\vartheta_{\infty}} Y^{\ubf_{\infty}}_{J_\infty}
$$
where each of the morphisms $\vartheta_n$ is a torsor under a
unipotent group of dimension $\dim(\Par_{J_n}) -
\dim(\Par_{J_{n-1}})$.

By~\Zastorsor\ the forgetful morphism $\pi\colon Z_{J_{\infty}}
\to \Par_{J_{\infty}} \times \Par_{K_{\infty}}$ is smooth and
surjective of relative dimension $\dim(G) -
2\dim(\Par_{J_{\infty}})$. The inverse image of $\SP^{u_{\infty}}
\subset \Par_{J_{\infty}} \times \Par_{K_{\infty}}$ under~$\pi$ is
nothing but $Y^{u_{\infty}}_{J_{\infty}} =
Y^{\ubf_{\infty}}_{J_{\infty}}$ as all pairs $(P,Q) \in
\SP^{u_{\infty}}$ are in good position.

In particular, we see that $Y^{\ubf_{\infty}}_J$ and hence $Y^{\ubf}_J$
are nonempty.

\secstart{Proposition}:\label{\Dimensionformula} {\sl
For $\ubf \in \Tscr(J)$ let $u_{\infty} \in {}^JW$ be the corresponding
element as in~\Combinatlemma. Then
$$
\codim(Y^{\ubf}_J, Z_J) = \dim(\Par_J) - \ell(u_{\infty})\, .
$$
}

\proof: By \Seedimension\ we have
$$\eqalign{
\dim(Y^{\ubf}_J) &= \dim(\Par_{J_{\infty}}) - \dim(\Par_J) +
\dim(Y^{\ubf_{\infty}}_{J_{\infty}}) \cr &=
\dim(\Par_{J_{\infty}}) - \dim(\Par_J) + \dim(G) -
2\dim(\Par_{J_{\infty}}) + \dim(\SP^{u_{\infty}}) \cr &= \dim(G) -
\dim(\Par_J) - \dim(\Par_{J_{\infty}}) + \ell(u_{\infty}) +
\dim(\Par_{J_{\infty}}) \cr &= \dim(G) - \dim(\Par_J) +
\ell(u_{\infty})\, .\cr }$$ On the other hand, \Zastorsor\ implies
that $\dim(Z_J) = \dim(G)$ which proves our claim.

\secstart{}\label{\uisuinfty} Suppose given an element $u \in
{}^JW^K$, rational over~$\Fq$, with the property that $J = {}^uK =
{}^{ux}\delta(J)$. The case we have in mind is when $J =
J_\infty$, $K = K_\infty$ and $u = u_\infty$ for some element
$\ubf \in \Tscr(J_0)$ as in~\usequences.

We fix a triple $(P_0,Q_0,L_0)$ consisting of a parabolic subgroup
$P_0 \subset G$ of type~$J$, a parabolic subgroup $Q_0 \subset G$
of type~$K$, and a subgroup $L_0 \subset G$ such that
$\relpos(P_0,Q_0) = u$ and such that $L_0$ is a common Levi
subgroup of $P_0$ and~$Q_0$. Such triples exists (rationally
over~$\Fq$) because $G$ is quasi-split and $J$, $K$ and~$u$ are
all defined over~$\Fq$. Note in particular that $F(P_0) = P_0$ and
$F(Q_0) = Q_0$.

Let $X^u_J$ be the $\Fq$-scheme whose $S$-valued points are the
elements $g \in G^1(S)$ such that
\assertionlist
\assertionitem $\relpos(Q_0,{}^gP_0) = x$;
\assertionitem ${}^gL_0 = L_0$;
\assertionitem $L_0$ is a Levi subgroup of~${}^gP_0$.

We have an action of $L_0$ on~$X^u_J$ by left multiplication. We claim
that this makes~$X^u_J$ an $L_0$-pseudo-torsor in the \'etale topology. (As
we will see below, $X^u_J$ is nonempty, so it is in fact a true
$L_0$-torsor.) To see this, suppose we have $g$, $h \in X^u_J(S)$. Then
$\relpos(Q_0,{}^gP_0) = x = \relpos(Q_0,{}^hP_0)$ and $L_0$ is a common
Levi subgroup of $Q_0$, ${}^gP_0$ and~${}^hP_0$. By~\Levigood\ this implies
that ${}^gP_0 = {}^h P_0$, hence the element $y := g^{-1}h$ lies
in~$P_0(S)$. But we also know that $y$ normalizes~$L_0$, so $y$ lies in the
normalizer of~$L_0$ inside~$P_0$, which is $L_0$ itself.

\secstart{} We have a second action of $L_0$ on~$X^u_J$, given on
points by $y \cdot g = ygF(y^{-1})$. (Note that $ygF(y^{-1})$ is again
in~$X^u_J$, as $F(y^{-1})$ is in $F(L_0) = L_0 \subset P_0$.) We denote
this action by~$\rho\colon L_0 \times X^u_J \to X^u_J$.

We have chosen $P_0$ and~$Q_0$ such that they are in good
position; in particular, $\Ref_{Q_0}(P_0) = P_0$. Hence if $g \in
X^u_J$ then $Q_0$ and ${}^gF\bigl(\Ref_{Q_0}(P_0)\bigr) =
{}^gF(P_0) = {}^gP_0$ are in standard position and we obtain a
well-defined morphism
$$
f\colon X^u_J \arr Y^u_J, \qquad g \asr [P_0,Q_0,g].
$$
Clearly $f$ is equivariant with respect to $L_0$-actions, where we take the
$\rho$-action on~$X^u_J$.

\secstart{Lemma}:\label{\connectYandX} {\sl
Notation and assumption as in~\uisuinfty. The morphism $G \times X^u_J \to
Y^u_J$ given on
points by
$$
(h,g) \mapsto \bigl({}^hP_0, {}^hQ_0, hgF(h^{-1})\bigr)
$$
is an epimorphism of fppf sheaves. In particular, if
$k$ is an algebraically closed extension field of~$\Fq$ then every
$G(k)$-orbit in~$Y^u_J(k)$ meets the image of~$X^u_J(k)$ under the
morphism~$f$.
}

\proof: The last assertion follows from the first because every
fppf covering of $\Spec(k)$ has a section. To prove the first
assertion, let $S$ be an $\Fq$-scheme and let $y \in Y^u_J(S)$.
After fppf-localization on~$S$ we may represent~$y$ by a triple
$(P,Q,g)$ in~$\tilde Y^u_J$. Possibly after a further localization
we can find an element $\gamma$ in~$G$ with ${}^\gamma Q = Q_0$.
Replacing $(P,Q,g)$ by $\bigl({}^\gamma P,{}^\gamma Q, \gamma
gF(\gamma^{-1})\big)$ we may from now on assume that $Q = Q_0$.

We know that $\relpos(P,Q_0) = u = \relpos(P_0,Q_0)$. Hence
fppf-locally on~$S$ we can find $\eta \in Q_0$ with ${}^\eta P =
P_0$. Replacing $(P,Q,g)$ by $\bigl({}^\eta P,{}^\eta Q, \eta
gF(\eta^{-1})\bigr)$ we arrive at the situation where $P = P_0$
and $Q= Q_0$.

The assumption that $\relpos\bigl(Q,{}^gF(P)\bigr) = x$ implies
that $Q$ and~${}^gF(P)$ are in good position, so fppf-locally
on~$S$ there is a common Levi~$M$ of $Q$ and~${}^gF(P) = {}^gP$.
There is a unique $v \in U_Q$ such that ${}^vM= L_0$. Replacing
$g$ by~$vg$ we get that $L_0$ is a common Levi of $P$, $Q$
and~${}^gP$. But then there is a unique $w \in U_P = U_{F(P)}$
with ${}^gL_0 = {}^w L_0$. Replacing $g$ by~$gw^{-1}$ we finally
arrive at a triple $(P,Q,g)$ that is in the image of~$X^u_J$
under~$f$.

\secstart{Lemma}:\label{\TransitiveonX} {\sl Notation and
assumption as in~\uisuinfty. The morphism $\Psi \colon L_0 \times
X^u_J \to X^u_J \times X^u_J$ given on points by $(y,g) \mapsto
\bigl(ygF(y^{-1}),g\bigr)$ is finite \'etale and surjective. In
particular, if $k$ is an separably closed field then $L_0(k)$
acts transitively on~$X^u_J(k)$. }

\proof: It follows from~\connectYandX\ and \Dimensionformula\ that
$X^u_J$ is nonempty. Choose a finite field extension $\Fq \subset
k$ such that $X^u_J(k) \neq \emptyset$. It suffices to show that
$\Psi$ is finite \'etale surjective after base change to~$k$. If $g
\in X^u_J(k)$ then we get an isomorphism $L_{0,k} \ariso
X^u_{J,k}$ by $z \mapsto zg$, and $\Psi$ becomes the morphism
$L_{0,k} \times L_{0,k} \to L_{0,k} \times L_{0,k}$ given by
$(y,z) \mapsto \bigl(yzgF(y^{-1})g^{-1},z\bigr)$.

Consider the morphism $h\colon L_{0,k} \to L_{0,k}$ given by $y
\mapsto yF(y^{-1})g^{-1}$. We claim that $h$ is finite \'etale and
surjective. In fact, it suffices to show this for the
morphism $h_1 \colon L_{0,k} \to L_{0,k}$ given by $y \mapsto
yF(y^{-1})$. By Lang's theorem, $h_1$ is surjective. The fibres
of~$h_1$ are principal homogeneous under (right multiplication
by)~$L_0(\Fq)$, and using [\Mats], Thm.~23.1 we find that $h_1$,
hence also~$h$, is finite faithfully flat. Looking at tangent spaces we see
that it is even \'etale.

We view $L_{0,k} \times L_{0,k}$ as a scheme over~$L_{0,k}$ via
the second projection. Note that $\Psi$ is a morphism
over~$L_{0,k}$. After base change over the morphism~$h$ we obtain
$$
h^\ast \Psi \colon L_{0,k} \times L_{0,k} \to L_{0,k} \times L_{0,k}
$$
given on points by $(c,d) \mapsto \bigl(cd
F((cd)^{-1})g^{-1},d\bigr)$. Writing $\mu\colon L_0 \times L_0 \to
L_0$ for the group law, we have an isomorphism $(\mu,\pr_2)\colon
L_0 \times L_0 \ariso L_0 \times L_0$, and we find that $h^\ast
\Psi = (h \times \id) \circ (\mu,\pr_2)$. Hence $h^\ast \Psi$ is
finite \'etale surjective, and since these properties are local for
the fppf topology, the lemma follows.

\secstart{Theorem}:\label{\GorbsinZJ} {\sl Let $\ubf \in \Tscr(J)$
and let $u_{\infty} \in {}^JW$ be the corresponding element as
in~\Combinatlemma. The $G$-scheme $Y^{\ubf}_J$ is equi-dimensional
of codimension $\dim(\Par_J) - \ell(u_{\infty})$ in~$Z_J$. The
group $G$ acts transitively on~$Y_J^\ubf$, in the sense that the
morphism
$$
G \times Y^{\ubf}_J \arr Y^{\ubf}_J \times Y^{\ubf}_J\, ,
\quad\hbox{given by}\quad
(g,y) \mapsto (y,g\cdot y)
$$
is an epimorphism of fppf sheaves.

In particular, for any algebraically closed extension $k$
of~$\F_q$ there is a natural bijection between the $G(k)$-orbits
in~$Z_J(k)$ and the set~${}^JW$. }

\proof: The dimension formula was proven in~\Dimensionformula. It
follows from \thetabijective\ that the $G$-action is transitive on
$Y^{\ubf}_J$ if and only if it is transitive on
$Y^{\ubf_{\infty}}_J = Y^{u_{\infty}}_J$. But this is the case
because of \connectYandX\ and~\TransitiveonX. The last assertion
now follows from~\Combinatlemma\ and~\ZJunionYJu.

\secstart{}\label{\ordinaryorbit} Note that in $Z_J$ there is a unique
orbit $Z_J^{\rm ord}$ of maximal dimension. Therefore this orbit is
open and dense in $Z_J$. We call $Z_J^{\rm ord}$ the {\it ordinary
orbit}. There is also a unique orbit of minimal dimension, which is 
therefore a closed orbit; it has dimension $\dim(G) - \dim(\Par_J)$.

In our applications to $F$-zips we shall have that $\delta = \id$ and 
$K = {}^{w_0}J$, where $w_0 \in W$ is the longest element. For $x \in 
{}^KW^J$ we take the minimal representative of $W_Kw_0W_J$; this element 
is in fact the unique element of maximal length in ${}^KW^J$. With these 
choices, if $R$ is an $\F_q$-algebra, $Z_J^{\rm ord}(R)$ consists of 
those points $[P,Q,g]$ in~$Z_J$ such that $\relpos(P,Q)$ is the element 
of maximal length $u_{\max} = x^{-1} \in {}^JW^K$. Indeed, if 
$\relpos(P,Q) = u_{\max}$ then $P$ and~$Q$ are in good 
position~\goodposition. Therefore $(P,Q) = (P_1,Q_1) = \cdots$, and we 
see that $u_{\infty} = u_0 = u_{\max}$. Conversely, under the bijection 
of Lemma~\Combinatlemma, the element $u_{\max} \in {}^JW^K \subset {}^JW$ 
corresponds to the constant sequence $\ubf = (u_{\max},u_{\max},\ldots)$.


\paragraph{Applications to F-zips}

\secstart{}\label{\transnot} Fix an integer $n \geq 0$. Let $V$ be an
$\Fp$-vector space of dimension~$n$. We shall apply the theory of Section~3
with $\Ghat = G := \GL(V)$. (So $q = p$.) Let $(W,I)$ be the Weyl
group with its set of simple reflections; see Example~\GLExa\ for
an explicit description.

If $S$ is an $\Fp$-scheme, write $V_S := V \otimes \Oscr_S$ and
$V_S^{(p)} := \Fr_S^\ast V_S = V_S \otimes_{\Oscr_S,\Fr_S}
\Oscr_S$. We have a canonical $\Oscr_S$-linear isomorphism $\xi_S
\colon V_S^{(p)} \arriso V_S$ by $(v \otimes x) \otimes y \mapsto
v \otimes x^py$ for $v \in V$ and $x$, $y$ local sections
of~$\Oscr_S$.

An $S$-valued point of~$G$ is given by an $\Oscr_S$-linear
automorphism $g$ of~$V_S$. The Frobenius endomorphism $F\colon G
\to G$ is given by $F(g) = \xi_S \circ g^{(p)} \circ \xi_S^{-1}$,
where $g^{(p)} = \Fr_S^\ast(g)$ is the automorphism $g \otimes
\id$ of~$V_S^{(p)}$.

We fix a function $\tau\colon \Z \to \Z_{\geq 0}$ with
$\sum_{i\in\Z} \tau(i) = n$. Let $C\udot$ be any filtration of
type~$\tau$ on~$V$. The stabilizer $\Stab(C\udot) \subset G$ is a
parabolic subgroup; its type $J \subset I$ is independent of the
choice of~$C\udot$. We refer to~$J$ as the {\it parabolic type
associated to~$\tau$\/}.

Let $w_0 \in W$ be the longest element, and set $K := {}^{w_0}J$.
Let $x \in {}^KW^J$ be the minimal representative of the double
coset $W_Kw_0W_J$. It is easily verified that $K = {}^xJ$, so we
are in the situation of~\JKx. (Note that $\delta$ is the
identity.)

\secstart{}\label{\GactsonX} Let $X_\tau$ be the scheme over~$\Fp$
whose $S$-valued points are the triples
$(C\udot,D\dodot,\phi\dodot)$ such that
$(V_S,C\udot,D\dodot,\phi\dodot)$ is an $F$-zip of type~$\tau$.
The group~$G$ acts on~$X_\tau$; on points:
$$
g \cdot (C\udot,D\dodot,\phi\dodot) =
\bigl(g(C\udot),g(D\dodot),\psi\dodot\bigr)\, ,
$$
where $\psi_i$ is the composition
$$
\eqalign{ \bigl(g(C^i)/g(C^{i+1})\bigr)^{(p)} \arriso
g^{(p)}\bigl((C^i/C^{i+1})^{(p)} \bigr) & \arrisoover{g^{(p),-1}}
\bigl(C^i/C^{i+1}\bigr)^{(p)}\cr & \arrisoover{\phi_i} D_i/D_{i-1}
\arrisoover{g} g(D_i)/g(D_{i-1})\, .}
$$

\secstart{Lemma}:\label{\XisoZ} {\sl With notation as above, there
is a $G$-equivariant isomorphism of $\Fp$-schemes $X_\tau \arriso
Z_J$.}

\proof: Consider the $\Fp$-scheme $\tilde{X}_\tau$ whose
$S$-valued points are the tuples
$$
\bigl(C\udot,\{A^i\}_{i\in\Z},D\dodot,\{B_i\}_{i\in\Z},\phi\dodot\bigr)
$$
with $(C\udot,D\dodot,\phi\dodot)$ in~$X_\tau$, with $\{A^i\}$ a
splitting of~$(C\udot)^{(p)}$ and $\{B_i\}$ a splitting
of~$D\dodot$. (By this we mean that $\{A^i\}_{i\in\Z}$ is a
collection of subspaces of~$V_S^{(p)}$ such that $(C^j)^{(p)} =
\oplus_{i \geq j} A^i$ for all~$j$; similarly for $\{B_i\}$.) We
have a forgetful morphism $\tilde{X}_\tau \to X_\tau$.

We may view~$X_\tau$, hence also~$\tilde{X}_\tau$, as schemes over
$\Par_J \times \Par_K$ by associating to
$(C\udot,D\dodot,\phi\dodot)$ the pair $(P,Q)$ with $P =
\Stab(C\udot)$ and $Q =\Stab(D\dodot)$. Let $\Uuni_{J,K}$, with
underlying scheme $F(\Uuni_J) \times \Uuni_K$, be the group scheme
over $\Par_J \times \Par_K$ as in~\defineZ. It acts from the left
on~$\tilde{X}_\tau$ over~$X_\tau$ by
$$
(u_1,u_2) \cdot
(C\udot,\{A^i\}_{i\in\Z},D\dodot,\{B_i\}_{i\in\Z},\phi\dodot) =
\bigl(C\udot,\bigl\{\xi_S^{-1}u_1^{-1}\xi_S(A^i)\bigr\}_{i\in\Z},
D\dodot,\bigl\{u_2(B_i)\bigr\}_{i\in\Z},\phi\dodot\bigr)\,
.
$$

The set of splittings of a filtration~$\Gamma\udot$ (descending or
ascending) is principal homogeneous under the unipotent radical of
the associated parabolic~$\Stab(\Gamma\udot)$. Using this fact it
readily follows that $X_\tau$ is the fppf quotient
of~$\tilde{X}_\tau$ modulo~$\Uuni_{J,K}$.

It remains to be shown that we have an isomorphism $\tilde{X}_\tau
\arriso \tilde{Z}_J$, equivariant with respect to both the
$G$-actions and the $\Uuni_{J,K}$-actions. Define $\alpha\colon
\tilde{X}_\tau \to \tilde{Z}_J$ by associating to an $S$-valued
point
$\bigl(C\udot,\{A^i\}_{i\in\Z},D\dodot,\{B_i\}_{i\in\Z},\phi\dodot\bigr)$
the triple $(P,Q,g)$ with $P = \Stab(C\udot)$ and $Q =
\Stab(D\dodot)$, and with $g \in G(S)$ the composition
$$
V_S \arrisoover{\xi_S^{-1}} V_S^{(p)} = \bigoplus_{i\in\Z} A^i \cong
\bigoplus_{i\in\Z} (\gr^i_C)^{(p)} \arrisoover{\phi\dodot}
\bigoplus_{i\in\Z} \gr_i^D \cong \bigoplus_{i\in\Z} B_i = V_S\, .
$$
By construction, $g\bigl(\xi_S(C\udot)^{(p)}\bigr)$ is in
opposition with~$D\dodot$; hence $\relpos\bigl(Q,{}^gF(P)\bigr) =
x$ and $\alpha$ is well-defined. It is straightforward to check
that $\alpha$ is equivariant with respect to the actions of~$G$
and~$\Uuni_{J,K}$.

Next we define a morphism $\beta\colon \tilde{Z}_J \to
\tilde{X}_\tau$. Start with an $S$-valued point $(P,Q,g) \in
\tilde{Z}_J$. Then $Q$ and~${}^gF(P)$ are two parabolics in
opposition, which means that $M := Q \cap {}^gF(P)$ is a common
Levi subgroup. Hence $L := {}^{g^{-1}}M$ is a Levi subgroup
of~$F(P)$. Now use the correspondences between parabolics and
flags, and between Levi subgroups and splittings of a flag. More
concretely, let $C\udot$ be the unique filtration of~$V_S$ of
type~$\tau$ such that $P = \Stab(C\udot)$, let $\{A^i\}$ be the
splitting of~$(C\udot)^{(p)}$ corresponding to the Levi subgroup
$L \subset F(P)$, let $D\dodot$ be the filtration of~$V_S$ of
type~$\tau$ corresponding to~$Q$, and let $\{B_i\}$ be the
splitting of~$D\dodot$ corresponding to the Levi subgroup~$M
\subset Q$. Because ${}^gL = M$, there is a permutation $\pi$
of~$\Z$ such that $g\bigl(\xi_S(A^i)\bigr) = B_{\pi(i)}$ for
all~$i\in\Z$. The assumption that $Q$ and~${}^gF(P)$ are in
opposition then implies that we in fact have
$g\bigl(\xi_S(A^i)\bigr) = B_i$ for all~$i$. Hence we can
define~$\phi_i$ to be the composition
$$
(\gr^i_C)^{(p)} \cong A^i \arrisoover{\xi_S \circ g} B_i \cong
D_i\, .
$$
Then
$\bigl(C\udot,\{A^i\}_{i\in\Z},D\dodot,\{B_i\}_{i\in\Z},\phi\dodot\bigr)$
is a well-defined element of~$\tilde{X}_\tau(S)$. As it is clear
from the construction that $\alpha$ and~$\beta$ are inverse to
each other, the lemma is proven.

\secstart{Theorem}: {\sl Let $k$ be an algebraically closed field
of characteristic $p>0$. Let $n \geq 0$ be an integer, let $G =
\GL_n$, and let $(W,I)$ be the Weyl group with its subset of
simple reflections. Let $\tau\colon \Z \to \Z_{\geq 0}$ be a
function with $\sum_{i\in\Z} \tau(i) = n$, and let $J \subset I$
be the associated parabolic type. Then there is a bijection
$$
\left\{\vcenter{
\setbox1=\hbox{{\rm $F$-zips of type $\tau$ over $k$}}
\hbox to\wd1{\hfill {\rm isomorphism classes of}\hfill}
\copy1
}\right\} \arriso {}^JW \cong W_J\backslash W\,
.\lformno\sublabel{\MainCorresp}
$$
In particular, every $F$-zip of type~$\tau$ is isomorphic to a
standard $F$-zip $\Mline_\tau^u \otimes_{\Fp} k$ as
in~\StandFzips, for a unique $u \in {}^JW$.}

\proof: The first statement is the conjunction of Thm.~\GorbsinZJ\
and the above lemma. For the second assertion one verifies that
$\Mline_\tau^u \otimes k$ corresponds, under~\MainCorresp,
precisely with the element $u \in {}^JW$.

\secstart{}\label{\OrdinaryFzip} We call the $\GL_n$-orbit of
$X_{\tau}$ corresponding via the bijection~\XisoZ\ to the open dense
orbit $Z_J^{\rm ord}$ in $Z_J$~\ordinaryorbit\ the {\it ordinary
orbit}. It parametrizes $F$-zips $(M,C\updot,D\dodot,\varphi\dodot)$
of type $\tau$ such that the filtration $C\updot$ and $D\dodot$ are in
opposition, i.e.,~the rank of $C^i \cap D_j$ is as small as
possible for all $i, j \in \Z$~\GLExa.

\secstart{}\label{\StandardGtriple} The standard $F$-zips
$\Mline_\tau^u$ defined in~\StandFzips\ correspond, under the
isomorphism of~\XisoZ, to certain standard triples $[P,Q,g]$
in~$Z_J$. As we shall discuss now, these can be defined
independent of the language of $F$-zips, for an arbitrary
reductive group~$\Ghat$ as in~\GhatSetup.

Let $L$ be the splitting field of~$G$. We choose an
$\F_q$-rational Borel pair $(T,B)$ of~$G$ such that $T$ is split
over~$L$. Via this choice we identify the Weyl group~$W$ with
$N_G(T)/T$. Moreover, we choose a set-theoretic section $s\colon
W(L) \arr N_G(T)\bigl(L\bigr)$. For $u \in {}^JW$ let $\ubf =
(u_0,u_1,\ldots) \in \Tscr(J)$ be the corresponding family
under~\Combinatlemma. We apply the definitions and the notation
of~\usequences ; in particular, $u = u_\infty$.

We denote by $(P^{\ubf}_{\infty},Q^{\ubf}_{\infty},g^{\ubf}) \in
\Ztilde_{J_\infty}(L)$ the triple satisfying: \indention{(b)}
\litem{(a)} $P^{\ubf}_{\infty} = {}^{u_\infty}P'$, where $P'$ is
the parabolic subgroup of type~$J_{\infty}$ containing~$B$;
\litem{(b)} $Q^{\ubf}_{\infty}$ is the parabolic subgroup of~$G_L$
of type $K_{\infty} = {}^x\delta(J_{\infty}) =
{}^{u_{\infty}^{-1}}J$ containing~$B$; \litem{(c)} $g^{\ubf} =
s\bigl((u_{\infty}x)^{-1}\bigr) \in N_G(T)\bigl(L\bigr)$.

Note that $K_{\infty} = {}^x\delta(J_{\infty})$ is already defined
over~$\F_q$; hence the same is true for $Q^{\ubf}_{\infty}$ and
$F(Q^{\ubf}_{\infty}) = Q^{\ubf}_{\infty}$. By definition we have
$$\eqalign{
\relpos\bigl(P^{\ubf}_{\infty}, F(Q^{\ubf}_{\infty})\bigr) &=
u_{\infty}, \cr \relpos(Q^{\ubf}_{\infty},
{}^{g^{\ubf}}P^{\ubf}_{\infty}) &= x\, . \cr }$$ Therefore,
$(P^{\ubf}_{\infty},Q^{\ubf}_{\infty},g^{\ubf}) \in
\Ytilde^{u_{\infty}}_J = \Ytilde^{\ubf_{\infty}}_J$.

Let $P^{\ubf}$ (resp.~$Q^{\ubf}$) be the unique parabolic of
type~$J$ (resp.~of type $K = {}^x\delta(J)$) containing
$P^{\ubf}_{\infty}$ (resp.~$Q^{\ubf}_{\infty}$).
Now~\thetasurjective\ implies that $(P^{\ubf},Q^{\ubf},g^{\ubf})
\in \Ytilde^{\ubf}_J$. We call the image
$[P^{\ubf},Q^{\ubf},g^{\ubf}] \in Y^{\ubf}_J$ the {\it standard
triple of type~$\ubf$ associated to $(T,B,s)$\/}. Another choice
of $(T,B,s)$ gives a point of~$Y^\ubf_J$ in the same $G$-orbit.

In the case $G = \GL_{n,\F_p}$, we have $L = \F_p$. Take $T$ to be
the diagonal torus, $B$ the Borel subgroup of upper triangular
matrices, and $s\colon W = S_n \to N_G(T)$ the map that associates
to a permutation the corresponding permutation matrix. Further,
fix a type~$\tau$ with $\sum \tau(i) = n$. The triple
$(P^\ubf,Q^\ubf,g^\ubf) \in \Ytilde_J$ then corresponds, under the
isomorphism as in the proof of~\XisoZ, to a point of the
scheme~$\tilde{X}_J$. It can be checked that this point is none
other than the standard $F$-zip~$\Mline^u_J$ together with the
obvious splitting of the filtrations $C\udot$ and~$D\dodot$ given
by the basis $\{e_1,\ldots,e_n\}$ of the underlying vector space.

\secstart{}\label{\somenot} Let $\tau\colon \Z \to \Z_{\geq 0}$ be a
function with finite support. Let $n := \sum_{i\in\Z} \tau(i)$. As
in~\transnot\ we fix an $\Fp$-vector space~$V$ of dimension~$n$, we set $G :=
\GL(V)$, and we let $(W,I)$ be the Weyl group of~$G$ with its set
of simple reflections. Let $J \subset I$ be the parabolic type
associated to~$\tau$.

Consider the scheme~$Z_J$. For $\ubf \in \Tscr(J)$ we have a
locally closed subscheme $Y_J^\ubf \hookrightarrow Z_J$, stable
under the action of~$G$, and the morphism
$$
\coprod_{\ubf \in \Tscr(J)} Y_J^\ubf \arr Z_J
\lformno\sublabel{\EOstratofZJ}
$$
is a bijective monomorphism.

The underlying reduced schemes $(Y_J^\ubf)_\red$ are irreducible and
non-singular, as $G$ acts transitively on them.

We claim that \EOstratofZJ\ is a stratification. More precisely, for
$\ubf$, $\vbf \in \Tscr(J)$, let us write $\vbf \preceq \ubf$ if
$Y_J^\vbf$ meets the Zariski closure of~$Y_J^\ubf$.
Then it follows from the general properties of orbits under a
group action (see e.g.\ [\TWDimStrat],~4.2) that ``$\preceq$'' is a
partial ordering on~$\Tscr(J)$ and that
$$
\overline{Y_J^\ubf} = \coprod_{\vbf \preceq \ubf} Y_J^\vbf\,
.\lformno\sublabel{\YJuClosure}
$$
This last identity has to be interpreted set-theoretically, or on
points with values in an algebraically closed field. As a slight
refinement, we shall prove in (4.11) below that if $\vbf \preceq \ubf$
then $Y_J^\vbf$ is in fact contained in the Zariski closure
of~$Y_J^\ubf$ as a subscheme.

\secstart{}\label{\EOstrat} 
{\it The Ekedahl-Oort stratification associated to an $F$-zip.\/} 
We retain the notation of \somenot. Let $\Mline$ be an $F$-zip of
type~$\tau$ over a  connected base scheme~$S$. Let $M$ be the
underlying locally free $\Oscr_S$-module, which  is of rank $n =
\dim(V)$. Let ${}^\#S \to S$ be the $G_S$-torsor of trivialisations
of~$M$; so if $T$ is a scheme over~$S$ then the $T$-valued points
of~${}\#S$ are the  isomorphisms $M_T \arriso V_T := (V \otimes_{\Fp}
\Oscr_T)$.

We have a canonical $G$-equivariant morphism
$$
\nu\colon {}^\#S \arr X_\tau \cong Z_J\, .
$$
As explained in the introduction, we think of this map as a ``mod~$p$
period map''.

For $u\in {}^JW$ corresponding to $\ubf \in \Tscr(J)$, define ${}^\#S^u := 
\nu^{-1}(Y_J^\ubf)$, which is a locally closed subscheme of~${}^\#S$,
preserved by the  action of~$G$. Now define $S^u \hookrightarrow S$ to
be the quotient of~${}^\#S^u$ by the  action of~$G$. Note that ${}^\#S
\to S$ is locally trivial for the Zariski topology, and  if $U \subset
S$ is an open subset over which we have a trivialisation $\alpha\colon
M_{|U} \arriso V_U$ then $S^u \cap U$ is just the pull-back of
${}^\#S^u \subset {}^\#S$  under the section $U \to {}^\#S$
corresponding to~$\alpha$. It is clear from the  construction that
$$
\coprod_{u \in {}^JW} S^u \arr S\lformno\sublabel{\EOstratofS}
$$
is a bijective monomorphism.

Let $\underline{\Mscr}$ be the universal $F$-zip over $X_{\tau} \cong Z_J$. The
locally closed subscheme $Y_J^\ubf$ is the locus where $\underline{\Mscr}$ is
fppf-locally isomorphic to the standard $F$-zip $\Mline^u_{\tau}$
using the notations \StandFzips\ with $q = p$. Hence it follows that
$S^u$ represents the subfunctor $\Sscr^u$ of $S$ which is defined by
the property that a morphism $g\colon T \ar S$ factors through
$\Sscr^u$ if and only if $g^*\Mline$ is fppf-locally isomorphic to
$\Mline^u_\tau \otimes_{\F_p} \Oscr_T$.

We refer to the subschemes $S^u \hookrightarrow S$ as the {\it
Ekedahl-Oort loci in~$S$ associated to the $F$-zip~$\Mline$\/} and
to~\EOstratofS\ as the {\it Ekedahl-Oort partition of~$S$
associated to~$\Mline$.} (We use the terms ``loci'' and
``partition'' because \EOstratofS\ is not, in general, a
stratification of~$S$. In fact, the closure of an irreducible
component of~$S^u$ need not be a union of components of EO-loci.)

\secstart{Definition}: In the above situation we say that the
$F$-zip~$\Mline$ is {\it isotrivial of type~$u$\/} if $S = S^u$.
We say that $\Mline$ is a {\it constant $F$-zip of type~$u$\/} if
it is isomorphic to $M_\tau^u \otimes_{\Fq} S$.

As a corollary of our method of proof we obtain the following
result.

\secstart{Corollary}:\label{\IsoTrivCor} {\sl Let $\Mline$ be an
$F$-zip of type~$\tau$ over~$S$. Then the following assertions are
equivalent:
\assertionlist
\assertionitem The $F$-zip $\Mline$ is
isotrivial of type~$u$.
\assertionitem There exists a faithfully
flat morphism $S' \to S$, locally of finite presentation, such
that $\Mline \otimes_S S'$ is a constant $F$-zip of type~$u$.

If $S$ is quasi-separated, these conditions are also equivalent
to:
\assertionitem Zariski-locally
on~$S$ there exists a faithfully flat quasi-finite morphism $S'
\to S$ of finite presentation such that $\Mline \otimes_S S'$ is a
constant $F$-zip of type~$u$.}

\proof: Our results in Section~3 show that if (1) holds then
$\Mline$ is fppf-locally constant; whence~(2). The equivalence of
(2) and~(3) in the quasi-separated case follows from [\EGA],
IV,~17.16.2.

\secstart{Remark}: Isotrivial $F$-zips are not, in
general, \'etale-locally constant. E.g., if the base scheme is the
spectrum of a field then in general we need a non-separable field
extension to trivialize the $F$-zip. As a concrete example, let $k$ 
be a field of characteristic~$p$, let $\gamma \in k$, and consider 
the $F$-zip $\Mline_\gamma$ with underlying module $M = k^5 = 
k\cdot e_1 + \cdots + k\cdot e_5$, with
$$
C^0 =M \;\supset\; C^1 = \Span(e_1,e_3) \;\supset\; C^2 = (0)
$$
and
$$
D_{-1} = (0) \;\subset\; D_0 = \Span(e_1,e_2,e_3) \;\subset\; D_1 = M\, ;
$$
with 
$$
\phi_0\colon (M/C^1)^{(p)} \arriso D_0
\quad\hbox{given by}\quad
\bar e_2^{(p)} \mapsto e_1\, ,\quad \bar e_4^{(p)} \mapsto e_2\, ,\quad 
\bar e_5^{(p)} \mapsto \gamma e_2 + e_3\, ,
$$
and with
$$
\phi_1\colon (C^1)^{(p)} \arriso M/D_0
\quad\hbox{given by}\quad
e_1^{(p)} \mapsto \bar e_4\, ,\quad e_3^{(p)} \mapsto \bar e_5\, .
$$
Then $\Mline_\gamma \cong \Mline_0$ over~$\kbar$, but to realize this 
isomorphism one has to extract a $p$th root of~$\gamma$.

\secstart{Lemma}: {\sl Let $\ubf$, $\vbf \in\Tscr(J)$ be elements
with $\vbf \preceq \ubf$. Let $\overline{Y_J^\ubf}$ denote the
scheme-theoretic Zariski closure of~$Y_J^\ubf$ (i.e., the
scheme-theoretic image of $Y_J^\ubf \to Z_J$). Then $Y_J^\vbf$
is contained in~$\overline{Y_J^\ubf}$ as subschemes of~$Z_J$.}

\proof: It suffices to show that if $R := k[t]/(t^n)$ with $k$ an
algebraically closed field, then any point $\mu\colon \Spec(R) \to
Y_J^\vbf$ factors through~$\overline{Y} := \overline{Y_J^\ubf}$.
But by~\IsoTrivCor, given an isotrivial $F$-zip of type~$\vbf$
over~$R$, then there is a faithfully flat extension $R \subset R'$
over which the $F$-zip is isomorphic to $M_\tau^\vbf \otimes_{\Fp}
R'$. In other words, if $m \colon \Spec(R') \to \Spec(\Fp) \to
Y_J^\vbf$ is the morphism corresponding to the constant $F$-zip
$M_\tau^\vbf \otimes R'$ then there is an element $g \in G(R')$
such that $\mu = g \cdot m$ in~$Z_J(R')$. Moreover,
by~\YJuClosure\ the point $\Spec(\Fp) \to Y_J^\vbf$ corresponding
to~$M_\tau^\vbf$ factors through~$\overline{Y}$, hence so does the
point~$m$. But $\overline{Y}$, as a closed subscheme of~$Z_J$, is
stable under the action of~$G$; hence $\mu \in \overline{Y}(R')$.


\paragraph{$F$-zips with additional structure}

\nnsection{}The purpose of this section is to discuss how the main
result of the previous section can be extended to $F$-zips with
additional structure. Ultimately one might wish to have a theory
of $F$-zips with $G$-structure, where $G$ is an arbitrary
reductive group. However, it is not clear to us how to define such
a notion in full generality. Therefore we restrict the discussion
to two simple examples.

\secstart{} Let $S$ be a scheme. Consider a pair $(M,\psi)$
consisting of a locally free $\Oscr_S$-module of finite rank,
together with a perfect pairing $\psi\colon M \otimes_{\Oscr_S} M
\to \Oscr_S$. Let $b_\psi\colon M \arriso M^\vee$ be the
isomorphism given on local sections by $m \mapsto \psi(-\otimes
m)$. For a locally direct summand $N \subset M$, we define
$N^\perp \subset M$ to be the kernel of the composite map $M
\arriso M^\vee \aerr N^\vee$, where the first map is~$b_\psi$. We
call $N$ isotropic if $N \subset N^\perp$; in that case $\psi$
induces a perfect pairing on~$N^\perp/N$. Note that
$N^{\perp\perp} = N$.

Now assume that either $\psi$ is symplectic, meaning that
$\psi(m,m) = 0$ for all local sections~$m$, or symmetric; we shall
consider the latter case only in characteristic $\neq 2$. A flag
$\Delta$ in~$M$ is called a symplectic (resp.\ orthogonal) flag if
for every $N \in \Delta$ we also have $N^\perp \in \Delta$. As
$\Delta$ is totally ordered, either $N$ or~$N^\perp$ is then
isotropic. We call a filtration symplectic (resp.\ orthogonal) if
the associated flag is.

Let $S$ be a scheme of characteristic~$p$. Consider a tuple
$\Mline = (M,\psi,C\udot,D\dodot,\phi\dodot)$ such that $\Mline'
:= (M,C\udot,D\dodot,\phi\dodot)$ is an $F$-zip over~$S$, with
$\psi$ a (perfect) symplectic or symmetric bilinear form on~$M$,
and such that the flags $C\udot$ and~$D\dodot$ are symplectic,
resp.\ orthogonal. Let $\tau$ be the type of~$C\udot$. Let $i \in
\Z$ be an index such that $\tau(i) \neq 0$. There is a unique
index $j\in\Z$ such that
$$
(C^i)^\perp = C^{j+1}
\quad\hbox{and}\quad
(C^{i+1})^\perp  = C^j\, ,
$$
and $b_\psi$ induces an isomorphism
$$
\alpha\colon \gr^j_C \arriso (\gr^i_C)^\vee\, .
$$
By an easy dimension count we then find that, for these same indices $i$
and~$j$, we have
$$
D_{i-1}^\perp = D_j
\quad\hbox{and}\quad
D_i^\perp = D_{j-1}\, ,
$$
and we get an isomorphism
$$
\beta\colon \gr^j_D \arriso (\gr_i^D)^\vee\, .
$$

\secstart{Definition}: Let $S$ be a scheme of characteristic~$p$,
with $p>2$ in the orthogonal case. By a {\it symplectic $F$-zip
over~$S$}, resp.\ an {\it orthogonal $F$-zip over~$S$}, we mean a
tuple $\Mline = (M,\psi,C\udot,D\dodot,\phi\dodot)$ as above, with
$\psi$ symplectic, resp.\ symmetric, such that for all indices $i$
and~$j$ as in the above discussion, the diagram
$$
\matrix{ (\gr^j_C)^{(p)} & \arrisoover{\phi_j} & \gr^j_D\cr
\addleft{\alpha^{(p)}} && \addright{\beta}\cr (\gr^i_C)^{(p),\vee}
& \allisounder{\phi_i^\vee} & (\gr_i^D)^\vee\cr}
$$
is commutative.

\secstart{}\label{\SymplOrth} Let $(V,\psi)$ be a finite
dimensional $\Fp$-vector space equipped with a perfect bilinear
pairing~$\psi$, assumed to be either symplectic or symmetric. If
$\psi$ is symmetric we assume that $p>2$ and also that $\dim(V)$
is {\it odd\/}.

In the symplectic case, set $G := \Sp(V,\psi)$; in the symmetric
case, $G := \SO(V,\psi)$. As usual, let $(W,I)$ be the Weyl group
with its set of simple reflections. We say that a type
$\tau\colon\Z \to \Z_{\geq 0}$ with support $i_1 < \cdots < i_r$
is {\it admissible\/} if $\tau(i_n) = \tau(i_{r+1-n})$ for
all~$n$. This is equivalent to the condition that for some
field~$k$ of characteristic~$p$, there exists a symplectic (resp.\
orthogonal) filtration~$C\udot$ of~$V_k$ of type~$\tau$. The
stabilizer $\Stab_G(C\udot)$ is then a parabolic subgroup
of~$G_k$; its type $J \subset I$ only depends on~$\tau$. We
call~$J$ the parabolic type associated to~$\tau$.

Define $X_\tau$ to be the $\Fp$-scheme whose $S$-valued points are
the triples $(C\udot,D\dodot,\phi\dodot)$ such that
$(V_S,\psi_S,C\udot,D\dodot,\phi\dodot)$ is a symplectic (resp.\
orthogonal) $F$-zip over~$S$. We let $G$ act on~$X_\tau$ by the
same rule as in~\GactsonX.

On the other hand, let $w_0$ be the element of maximal length in the
Weyl group $W$ of $G$ and let $x \in W_{w_0(J)}w_0W_J$ be the element
of minimal length. Consider the scheme~$Z_J$ defined in~\defineZ\
associated to $G$, $J$ and $x$. We claim that we again have a
$G$-equivariant isomorphism of $\Fp$-schemes $X_\tau \arriso Z_J$.
The proof of this is essentially the same as that of~\XisoZ,
provided we consider symplectic (resp.\ orthogonal) splittings of
the filtrations $(C\udot)^{(p)}$ and~$D\dodot$. We leave the
details to the reader. Note, however, that it is essential to have
a bijective correspondence between symplectic (resp.\ orthogonal)
flags and parabolic subgroups of~$G$, as well as a correspondence
between the symplectic (resp.\ orthogonal) splittings of a flag
and the Levi subgroups of the corresponding parabolic. Such a
correspondence fails for orthogonal groups in an even number of
variables, which is why we assume that in the orthogonal case,
$\dim(V)$ is odd.

\secstart{Corollary}:\label{\MainThmPol} {\sl Let $k$ be an
algebraically closed field of characteristic~$p$.

{\rm (\romannumeral1)} Let $G = \Sp(V,\psi)$ and $(W,I)$ be as
above; symplectic case. Let $\tau$ be an admissible type with
$\sum_{i\in\Z} \tau(i) = \dim(V)$ and with associated parabolic
type~$J \subset I$. Then there is a bijection
$$
\left\{\vcenter{\setbox0=\hbox{{\rm isomorphism classes of symplectic}}
\copy0
\hbox to\wd0{\hfill {\rm $F$-zips of type $\tau$ over $k$}\hfill}}\right\}
\arriso {}^JW \cong W_J\backslash W\, .
$$

{\rm (\romannumeral2)} Let $G = \SO(V,\psi)$ and $(W,I)$ be as
above; orthogonal case, with $\dim(V)$ odd. Let $\tau$ be an
admissible type with $\sum_{i\in\Z} \tau(i) = \dim(V)$ and with
associated parabolic type~$J \subset I$. Then there is a bijection
$$
\left\{\vcenter{\setbox0=\hbox{{\rm isomorphism classes of orthogonal}}
\copy0
\hbox to\wd0{\hfill {\rm $F$-zips of type $\tau$ over $k$}\hfill}}\right\}
\arriso {}^JW \cong W_J\backslash W\, .
$$}

Note that in this result the $\Fp$-structure on~$G$ plays no role,
as $(W,I)$ only depends on~$G_k$, which in turn only depends on
$\dim(V)$.

\secstart{Remark}: As remarked at the beginning of this section,
it is not clear to us how to define the notion of an $F$-zip with
$G$-structure, for $G$ an arbitrary reductive group. It is
possible, though, to obtain rather complete results for $F$-zips
equipped with an action of a semi-simple algebra and a hermitian
form. For Dieudonn\'e modules this was carried out in~[\BMGSAS].

\secstart{}\label{\SpSOinGL} Slightly changing notation, let $G_1
:= \Sp(V,\psi)$, resp.\ $G_1 := \SO(V,\psi)$ be the reductive
group over~$\Fp$ considered in~\SymplOrth. Let $G_2 := \GL(V)$.
Let $(W_i,I_i)$ be the Weyl group of~$G_i$.

If $\Par(G_i)$ is the scheme of parabolic subgroups of~$G_i$ then
we have a canonical morphism $\Par(G_1) \hookrightarrow
\Par(G_2)$; in terms of symplectic (resp.\ orthogonal)
flags~$\Delta$ in~$V$ it sends $\Stab_{G_1}(\Delta)$ to
$\Stab_{G_2}(\Delta)$. As $W_i$ can be identified with the set of
$G_i$-orbits in $\Par(G_i)_\emptyset^2$ (over any separably closed
field), we obtain a natural homomorphism $\iota\colon W_1 \to
W_2$, which is in fact injective.

Let $k = \kbar$. Let $\Mline$ be a symplectic (resp.\ orthogonal)
$F$-zip over~$k$ with $\dim(M) = \dim(V)$. Write $\Mline'$ for the
underlying $F$-zip, obtained by forgetting the form~$\psi$. Let
$J_i \subset I_i$ be the parabolic type associated to the
type~$\tau$ in the group~$G_i$. Then $\iota$ maps ${}^{J_1}W_1$
into ${}^{J_2}W_2$. If $u_1 \in {}^{J_1}W_1$ is the element
corresponding to~$\Mline$ under~\MainThmPol, and $u_2 \in
{}^{J_2}W_2$ is the element corresponding to~$\Mline'$
under~\MainCorresp\ then we have the relation $\iota(u_1) = u_2$.

For a more precise statement, consider the schemes $Z_{J_1}^{(1)}
\arriso X_\tau^{(1)}$ formed with respect to the group~$G_1$ and
the subset $J_1 \subset I_1$ and the schemes $Z_{J_2}^{(2)}
\arriso X_\tau^{(2)}$ formed with respect to~$G_2$ and $J_2
\subset I_2$. Then the forgetful morphism $\Mline \mapsto \Mline'$
defines a closed immersion $\alpha \colon X_\tau^{(1)}
\hookrightarrow X_\tau^{(2)}$. If $u_1 \in {}^{J_1}W_1$
corresponds to the sequence $\ubf_1 \in \Tscr(J_1)$ and $u_2 :=
\iota(u_1)$ corresponds to $\ubf_2 \in \Tscr(J_2)$ then it can be
shown that $\alpha$ induces an isomorphism between the subscheme
$Y_{J_1}^{\ubf_1} \hookrightarrow Z_{J_1}^{(1)} \cong
X_\tau^{(1)}$ and the subscheme $Y_{J_2}^{\ubf_2} \cap
Z_{J_1}^{(1)} \hookrightarrow Z_{J_2}^{(2)} \cong X_\tau^{(2)}$.



\paragraph{F-zips coming from geometry}

\secstart{} Let $f\colon X \to S$ be a morphism of schemes in
characteristic~$p$. We denote by $\Frob_S\colon S \to S$ the
absolute Frobenius. By definition of the relative
Frobenius~$F_{X/S}$ we have a commutative diagram
$$\matrix{
X & \arrover{F_{X/S}} & X^{(p)} & \arrover{\sigma} & X \cr &
\kern-10pt {\scriptstyle f}\searrow & \addright{f^{(p)}} & &
\addright{f} \cr & & S & \arrover{\Frob_S} & S \cr }$$ where the
square is cartesian.

Now assume that $f$ is smooth. Recall that we have two spectral
sequences converging to the de Rham cohomology $H\updot_{\DR}(X/S)
= {\bf R}\updot_*(\Omega\updot_{X/S})$, namely the Hodge-de Rham
spectral sequence
$${}_HE^{ab}_1 = R^bf_*(\Omega_{X/S}^a) \implies H^{a+b}_{\DR}(X/S)$$
and the conjugate spectral sequence
$${}_{\rm conj}E^{ab}_2 = R^af_*\bigl(\Hscr^b(\Omega\updot_{X/S})\bigr)
\implies H^{a+b}_{\DR}(X/S)\, .$$

Moreover, there is a unique isomorphism of graded
$\Oscr_{X^{(p)}}$-modules
$$
\Cscr^{-1}\colon \bigoplus_{i\geq 0}\Omega^i_{X\upp/S} \arriso
\bigoplus_{i\geq 0}\Hscr^i\bigl(F_*(\Omega\updot_{X/S})\bigr)\, ,
\lformno\sublabel{\CartIso}
$$
the (inverse) Cartier isomorphism, which satisfies
$$
\eqalign{
\Cscr^{-1}(1) &= 1 \cr
\Cscr^{-1}(d\sigma^{-1}(x)) &= \hbox{class of $x^{p-1}dx$} \cr
\Cscr^{-1}(\omega \wedge \omega') &= \Cscr^{-1}(\omega) \wedge
\Cscr^{-1}(\omega')\, . \cr
}$$

\secstart{} Let $f\colon X \to S$ be a smooth and proper morphism.
We say that~$f$ {\it satifies condition {\rm (D)}\/} if the following
two conditions hold: \indention{(b)} \litem{(a)} The $\Oscr_S$-modules
$R^bf_*(\Omega^a_{X/S})$ are locally free of finite rank for all
$a,b \geq 0$. \litem{(b)} The Hodge-de Rham spectral sequence
degenerates at~$E_1$.

If $f$ is satisfies (D), the formation of the Hodge-de Rham spectral
sequences commutes with base change $S' \to S$.

\secstart{}\label{\conjugateSS} Let $f\colon X \to S$ be a smooth
morphism of schemes of characteristic~$p$. For $a$, $b \in
\Z_{\geq 0}$ the (inverse) Cartier isomorphism $\Cscr^{-1}$
of~\CartIso\ defines an isomorphism
$$
R^af\upp_*(\Omega_{X\upp/S}^b) \arriso {}_{\rm conj}E^{ab}_2 =
R^af_*\bigl(\Hscr^b(\Omega\updot_{X/S})\bigr)\, .
$$
If further the $\Oscr_S$-modules $R^pf_*(\Omega^q_{X/S})$ are flat (e.g.\
if $f$ satisfies condition~(D)), we get an isomorphism
$$
\varphi^{ab}\colon \Frob_S^*R^af_*(\Omega^b_{X/S}) =
\Frob_S^*({}_HE^{ba}) \arriso {}_{\rm conj}E^{ab}_2 =
R^af_*\bigl(\Hscr^b(\Omega\updot_{X/S})\bigr)\, .
$$
Using this, one can show (e.g. [\KatzDiffEq],~2.3.2), that if $f$ is
satisfies condition (D), the conjugate spectral sequence degenerates at~$E_2$
and that its formation commutes with arbitrary base change.

\secstart{} We list some examples of morphisms that satisfy condition
(D). As usual, $S$ is a scheme of characteristic~$p$. 
\assertionlist
\assertionitem Any abelian scheme $f\colon A \to S$ is
satisfies (D). (Degeneracy of the Hodge-de Rham spectral sequence at~$E_1$
can be proven as in~[\Oda], Prop.~5.1.)
\assertionitem Any smooth proper curve $f\colon C \to S$ satisfies
(D). (Use the previous example.) 
\assertionitem Any K3-surface $X \to S$ satisfies (D). (This follows
from [\DelRelKKK], Prop.~2.2.) 
\assertionitem Every smooth complete intersection in
the projective space $\P^n_S$ satisfies (D) as a scheme over~$S$.
(See [\DelCohInt], 
Thm.~1.5.)
\assertionitem Let $f\colon X \to S$ be a smooth proper morphism
such that $(F_{X/S})_*(\Omega^{\bullet}_{X/S})$ is decomposable
(i.e., isomorphic in the derived category to a complex with zero
differential). Then $f$ satisfies (D) by results of Deligne and Illusie,
see [\DelIll], Cor.~4.1.5. Moreover,
this condition is satisfied if $\dim(X/S) < p$ and $f$ admits a
smooth lifting $\tilde{f}\colon \Xtilde \to \Stilde$ with
$\Stilde$ a flat $\Z/p^2\Z$-scheme (loc.\ cit.,~3.7).

\secstart{}\label{\HodgeWitttoFzip} Let $f\colon X \to S$ be a
morphism satisfying (D). Fix an integer $n$ with $0 \leq n \leq
2\dim(X/S)$. We associate to~$f$ an $F$-zip $(M,C\updot,D\dodot,
\varphi\dodot)$ over~$S$ as follows: Set $M = H^n_{\DR}(X/S)$. Let
$C\updot$ be the Hodge filtration on~$M$, and define the
filtration~$D\dodot$ by $D_i = {}_{\rm
conj}F^{n-i}H^n_{\DR}(X/S)$. Finally, let
$$
\varphi_i := \varphi^{n-i,i}\colon (\gr^i_C)^{(p)} =
\Frob_S^*R^{n-i}f_*(\Omega^i_{X/S}) \arriso  \gr_i^D =
R^{n-i}f_*(\Hscr^i(\Omega\updot_{X/S}))\, ,
$$
where $\varphi^{n-i,i}$ is the isomorphism defined in~\conjugateSS.

Note that $C\updot$ and~$D\dodot$ are filtrations in the sense
of~\defineFil. This follows from the fact that both the Hodge-de
Rham spectral sequence and the conjugate spectral sequence are
compatible with base change. (Use the fact that a homomorphism
$\iota\colon N \to M$ of $\Oscr_S$-modules with $M$ locally free
of finite type, makes~$N$ into a direct summand of~$M$ if and only
if $\iota$ stays injective after arbitrary base change.)

We obtain a functor $\FZ(n)$ from the category of $S$-schemes
$f\colon X \to S$ that satisfy condition (D) into the category of
$F$-zips over~$S$. This functor is compatible with base change $S'
\to S$.

\secstart{} In the situation of~\HodgeWitttoFzip, if $S$ is smooth over
some scheme~$T$ then we have an additional structure on the $F$-zip
$H^n_\DR(X/S)$, viz.\ a Gauss-Manin connection~$\nabla$ relative
to~$T$. The filtration~$D\dodot$ is horizontal with respect
to~$\nabla$, the filtration~$C\updot$ satisfies Griffiths
transversality, and the maps~$\varphi_i$ are horizontal with respect
to the canonical connection on $(\gr^i_C)^{(p)}$ and the connection
induced by~$\nabla$ on~$\gr_i^D$. In this paper we shall make no
attempt to further exploit this structure.

\secstart{} There is also a logarithmic variant. For this we use
the language of logarithmic schemes, as for instance in Kato's 
paper~[\KatoLog]. Let $f\colon
(X,\Mscr) \to (S,\Nscr)$ be a morphism of schemes with fine
log-structures in characteristic~$p$ and denote by
$\omega_{X/S}^{\bullet}$ the logarithmic de Rham complex. As in
the non-logarithmic case, there are two spectral sequences
converging to the logarithmic de Rham cohomology
$H\updot_{\DR,\log}(X/S) = {\bf R}\updot f_*(\omega\updot_{X/S})$:
$$\eqalign{
{}_HE^{ab}_1 = R^bf_*(\omega_{X/S}^a) &\implies
H^{a+b}_{\DR,\log}(X/S),\cr {}_{\rm conj}E^{ab}_2 =
R^af_*\bigl(\Hscr^b(\omega\updot_{X/S})\bigr) &\implies
H^{a+b}_{\DR,\log}(X/S)\, .}$$

If $f$ is log-smooth and of Cartier type, there exists also a
logarithmic variant $\Cscr^{-1}$ of the Cartier isomorphism.

Similarly as above, we say that $f$ {\it satisfies condition {\rm (D)}
\/} if the following conditions are satisfied:
\indention{(b)}
\litem{(a)} The log-structures $\Mscr$ and~$\Nscr$ are fine and the
morphism $f$ is log-smooth and of Cartier type. Its
underlying scheme morphism is proper.
\litem{(b)} The logarithmic Hodge-de Rham spectral sequence degenerates
at level~$E_1$.
\litem{(c)} The $\Oscr_S$-modules $R^bf_*\omega^a_{X/S}$ are locally free.

If the log-structures $\Mscr$ and~$\Nscr$ are trivial (or more
general if $f$ is a strict morphism of log-schemes) then $f$ satisfies
condition (D) if and only if the underlying scheme morphism is
satisfies (D). A nontrivial example for a morphism of
log-schemes satisfying condition (D) is the following case: Let $S$ be
the spectrum of a discrete valuation ring and let $X$ be a complete
intersection in a projective space over~$S$. Assume that $X$ is a regular, flat
over~$S$, and that its special fibre is a divisor with normal
crossings. Then the structure morphism $f\colon X \to S$ satisfies
condition~(D) if we endow $X$ and~$S$ with their natural log-structures.

Again one can show that condition~(b) and the existence of the
Cartier isomorphism imply that the conjugate spectral sequence
degenerates at~$E_2$. Moreover, condition~(3) and the existence of
the Cartier isomorphism then imply that the formation of the
logarithmic Hodge-de Rham spectral sequence and of the
logarithmic conjugate spectral sequence commute with arbitrary
base change.

For a log-smooth morphism $f\colon (X,\Mscr) \arr (S,\Nscr)$ of
fine log-schemes, the sheaf of logarithmic differentials
$\omega^1_{X/S}$ is locally free of finite type. If its rank is
constant we call this rank the relative dimension of $(X,\Mscr)$
over $(S,\Nscr)$ and denote it by~$\dim(X/S)$. (Note that in
general the underlying scheme morphism of~$f$ need not even be
flat.) If $f$ now satisfies condition (D) then we obtain, as in the
non-logarithmic case, an $F$-zip structure on
$H^n_{\DR,\log}(X/S)$ for every integer~$n$ with $0 \leq n \leq
2\dim(X/S)$.

\secstart{}\label{\Ordinarysetup} Let $f\colon X \to S$ be a
morphism satisfying (D). Fix $0 \leq n \leq 2\dim(X/S)$, and denote by
$\FZ(n)(f) = (M,C\updot,D\dodot,\varphi\dodot)$ the corresponding
$F$-zip with $M = H^n_\DR(X/S)$. We assume that $N(n) = \rk_{\Oscr_S}(M)$ is
constant on~$S$. Let $J(n)$ be the parabolic type associated to~$C\updot$.
Let $W(n) = S_{N(n)}$ be the Weyl group of $\GL_{N(n)}$ and let
$w(n)_{\max}$ be the unique maximal element in ${}^{J(n)}W(n)$ with
respect to the Bruhat order.

The Ekedahl-Oort locus $S^{w(n)_{\max}}$ corresponding to
$w(n)_{\max}$ and the choice of~$n$ is an open subscheme of~$S$. We
set
$$
S_{\ord} = \bigcap_n S^{w(n)_{\max}}\, .
$$

\secstart{Proposition}: {\sl In the situation of~\Ordinarysetup\
we have $S_{\ord} = S$ if and only if for every geometric point
$\sbar$ of~$S$ the $\kappa(\sbar)$-scheme $X_{\sbar}$ is ordinary
in the sense of [\IllRay],~4.12. }

\proof: As $S_{\ord} \subset S$ is open, we can assume that $S =
\Spec(k)$ for an algebraically closed field~$k$. By [\IllRay],~4.13,
$X$ is ordinary if and only if Hodge filtration and conjugate
filtration in $H^n_{\DR}(X/k)$ are in opposition, i.e., if and
only if their relative position is equal to the maximal element in
${}^{J_n}W_n^{K_n}$ where $K_n = w_0(J_n)$ and where $w_0$ is the
maximal element in~$W_n$. By~\OrdinaryFzip\ this is the case if and
only if the isomorphism type of $\FZ_n(f)$ corresponds,
via~\MainCorresp, to the element~$w_{n,\max}$.

\secstart{} Let $X$ be a (log-)smooth projective variety over an
algebraically closed field~$k$ such that $X \to \Spec(k)$ satisfies
condition (D). As suggested by the title of this paper, we may think
of the $F$-zip structure on the de Rham cohomology as a discrete
invariant of~$X$. As such, this contains certain discrete
invariants previously studied by other authors, such as the 
$a$-number defined
by van der Geer and Katsura in~[\GeerKats]. 
More precisely, if $u_\infty \in {}^JW$ is the element classifying 
the $F$-zip $H^m_\DR(X/k)$, and if $\ubf = (u_0,u_1,\ldots)$ is 
the sequence corresponding to~$u_\infty$ via~\Combinatlemma, then
the $a$-number only depends on~$u_0$, which is the relative 
position of the Hodge and the conjugate filtration.

\secstart{} {\it $F$-zips and Shimura varieties of PEL-type.\/}
Let $\Dscr = \bigl(B,\star,V,\lrangle, O_B, \Lambda, h\bigr)$
denote a Shimura-PEL-datum, integral and unramified at a
prime~$p$, let $\Gbf$ its associated reductive group over~$\Q$,
and $[\mu]$ denotes the associated conjugacy class of cocharacters
of~$\Gbf$. By this we mean that 
\indention{$\bullet$}
\litem{$\bullet$} $B$ is a finite-dimensional semi-simple
$\Q$-algebra, such that $B_{\Qp}$ is isomorphic to a product of
matrix algebras over unramified extensions of~$\Qp$;
\litem{$\bullet$} $\star$ is a $\Q$-linear positive involution
on~$B$; 
\litem{$\bullet$} $V \not= 0$ is a finitely generated left
$B$-module; 
\litem{$\bullet$} $\lrangle$ is a nondegenerate
alternating $\Q$-valued form on~$V$ such that $\langle bv,w
\rangle = \langle v,b^*w \rangle$ for all $v,w \in V$ and $b \in
B$; 
\litem{$\bullet$} $O_B$ is a $\star$-invariant
$\Z_{(p)}$-order of~$B$ such that $O_B \otimes \Zp$ is a maximal
order of $B \otimes \Qp$; 
\litem{$\bullet$} $\Lambda$ is an
$O_B$-invariant $\Zp$-lattice in~$V_{\Qp}$, such that
$\lrangle\restricted{\Lambda \times \Lambda}$ is a perfect pairing
of $\Zp$-modules; 
\litem{$\bullet$} $\Gbf$ is the $\Q$-group of
$B$-linear symplectic similitudes of~$V$, i.e., for any
$\Q$-algebra~$R$ we have
$$
\Gbf(R) = \bigl\{g \in \GL_B(V \otimes R) \bigm| \langle gv,gw 
\rangle = c(g) \cdot \langle
v,w \rangle\ \hbox{for some}\ c(g) \in R\cross \bigr\}\, ;
$$
\litem{$\bullet$} $h\colon \Res_{\C/\R}(\G_{m,\C}) \ar \Gbf_{\R}$
is a homomorphism defining a complex structure on $V_{\R}$ which
is compatible with~$\lrangle$; 
\litem{$\bullet$} $[\mu]$ is the
$\Gbf(\C)$-conjugacy class of the cocharacter~$\mu_h$ associated
to~$h$ (cf.\ [\DelVarSh],~1.1.1). Then $V_{\C}$ has only weights 0 and~1
with respect to any $\mu \in [\mu]$.

We assume that $p > 2$ if $\Gbf$ is not connected.

Let $E$ be the associated reflex field, i.e., the field of
definition of~$[\mu]$. It is a finite extension of~$\Q$. Fix an
embedding of the algebraic closure~$\Qdbar$ of $\Q$ in~$\C$ into
some algebraic closure $\Qpbar$ of~$\Qp$. Via this embedding we
can consider~$[\mu]$ as a $\Gbf(\Qpbar)$-conjugacy class of
cocharacters. Denote by $v\vert p$ the place of~$E$ given by the
chosen embedding $\Qdbar \air \Qpbar$ and write $E_v$ for the
$v$-adic completion of~$E$. Let $\kappa = \kappa(v)$ be its
residue class field.

Further fix an open compact subgroup $K^p \subset \Gbf(\Abf^p_f)$
and denote by ${\tt A}_{\Dscr,K^p}$ the associated moduli space,
defined by Kottwitz in~[\KottPts]. We assume that $K^p$ is sufficiently
small such that ${\tt A}_{\Dscr,K^p}$ is representable. It is then
a smooth equi-dimensional quasi-projective scheme over the
localization of $O_E$ in~$p$. It classifies tuples
$(A,\lgbar,\iota,\hgbar)$ where \indention{$\bullet$}
\litem{$\bullet$} $A$ is an abelian scheme up to
prime-to-$p$-isogeny; \litem{$\bullet$} $\lgbar$ is a
$\Q$-homogeneous polarization of~$A$ containing a polarization
$\lambda \in \lgbar$ of degree prime to~$p$; \litem{$\bullet$}
$\iota\colon O_B \ar \End(A) \otimes_{\Z} \Z_{(p)}$ is an
involution preserving $\Z_{(p)}$-algebra homomorphism where the
involution is $\star$ on~$O_B$ and the Rosati involution given
by~$\lgbar$ on $\End(A) \otimes_{\Z} \Z_{(p)}$; \litem{$\bullet$}
$\hgbar$ is a $K^p$-level structure.

\noindent Further $(A, \lgbar, \iota, \hgbar)$ should satisfy a
determinant condition; see [\KottPts],~\pz5 or [\RapZink],~3.23~a) for a
precise formulation. We denote by ${\tt A}_0$ the reduction ${\tt
A}_{\Dscr,K^p} \otimes \kappa$ at~$v$.

\secstart{} We denote by $\Ghat$ the reductive $\Fp$-group of
$O_B/pO_B$-linear symplectic similitudes of $\Lambda_0 :=
\Lambda/p\Lambda$ and let $G$ be its identity component. Via the
canonical bijection of $\Gbf(\Qdbar_p)$-conjugacy classes of
cocharacters and $G(\Fdbar_p)$-conjugacy classes of cocharacters
we consider $[\mu]$ as a $G(\Fdbar_p)$-conjugacy class of
cocharacters. Its field of definition is~$\kappa$. Let $(W,I)$ be
the Weyl group of~$G$ together with its set of simple reflections,
and let $J \subset I$ be the subset of simple reflections
corresponding to~$[\mu]$. Then $J$ is defined over~$\kappa$.

\secstart{} Let $S$ be a $\kappa$-scheme and let
$(A,\lgbar,\iota,\hgbar)$ be an $S$-valued point of~$\Att_0$.
Every abelian scheme satisfies condition (D). We set $M =
H^1_{\DR}(A/S)$. By~\HodgeWitttoFzip\ we obtain the structure of
an $F$-zip on~$M$. The filtration $C\updot$ (resp.~$D\dodot$) is
of the form $M = C^0 \supset C^1 \supset C^2 = (0)$ (resp.\ $(0) =
D_{-1} \subset D_0 \subset D_1 = M$), where $C^1$ and~$D_0$ are
locally direct summands of rank equal to $\dim(A/S)$. Moreover,
the submodules $D_0$ and~$C^1$ are $O_B/pO_B$-invariant and
totally isotropic with respect to the perfect alternating form
induced by any $\lambda \in \lgbar$ which is of order prime
to~$p$.

\secstart{Lemma}:\label{\DRloctrivial} {\sl Locally for the
\'etale topology the two skew Hermitian modules with
$O_B/pO_B$-action $M$ and $\Lambda_{0,S} = \Lambda_0 \otimes{\Fp}
\Oscr_S$ are isomorphic.}

\proof: This is a special case of~[\RapZink],~3.16.

\secstart{} We define two smooth coverings ${}^{\#}{\tt A}_0$ and
$\tilde{\tt A}_0$ of~${\tt A}_0$ as follows: For every
$\kappa$-scheme~$S$ the $S$-valued points of ${}^{\#}{\tt A}_0$
are given by tuples $(A,\lgbar,\iota,\hgbar,\alpha)$ where
$(A,\lgbar,\iota,\hgbar) \in {\tt A}_0(S)$ and where $\alpha$ is
an $O_B/pO_B$-linear symplectic similitude $H^1_{\DR}(A/S) \arriso
\Lambda_{0,S}$.

The $S$-valued points of $\tilde{\tt A}_0$ are given by tuples
$(A,\lgbar,\iota,\hgbar,\alpha, C',D')$ with
$(A,\lgbar,\iota,\hgbar,\alpha) \in {}^{\#}{\tt A}_0$ and where
$C'$ and~$D'$ are $O_B/pO_B$-invariant totally isotropic
complements of $C^1$ and~$D_0$, respectively.

It follows from~\DRloctrivial\ that ${}^{\#}{\tt A}_0$ is a torsor
for the \'etale topology over~${\tt A}_0$ under the smooth group
scheme~$G$. Furthermore, because Zariski-locally on~$S$ we can
always find complements $C'$ and~$D'$ as above, $\tilde{\tt A}_0$
is a torsor over ${}^{\#}{\tt A}_0$ under the smooth unipotent
group scheme $\Uuni_{J,K}$ defined in~\defineZ, where $J$ and~$K$
are the parabolic types of the filtrations $C\updot$
and~$D\dodot$, respectively.

We relate this to $\Ztilde_J$ as defined in~\DefineZtilde. The
scheme $\Ztilde_J$ depends on some automorphism $\delta$ of the
Weyl group of~$G$ which takes into account that our group $\Ghat$
might be disconnected. Hence we will write $\Ztilde_{J,\delta}$.
Moreover we let $\Ztilde'_J$ be the disjoint union of the schemes
$\Ztilde_{J,\delta}$ for the various possible~$\delta$. We can do
this also for the schemes $Z_J = Z_{J,\delta}$ and obtain a
scheme~$Z'_J$.

For every $S$-valued point $(A,\dots)$ of $\tilde{\tt A}_0$ we obtain
an $F$-zip with underlying $\Oscr_S$-module $H^1_\DR(A/S)$ with
additional structures and with splittings for their filtrations. As in
the proof of~\XisoZ\ we can associate to this $F$-zip an $S$-valued
point of~$\Ztilde'_J$. By passing to the quotients, we obtain a morphism
$$
\pi\colon \Att_0 \arr \bigl[G\backslash Z'_J\bigr]
$$
where on the right hand side we have the quotient stack.

Note that it follows from [\BMSTPEL],~4.1, that we can decompose
${\Att}_0$ into the special fibres of individual Shimura varieties
such that $\pi$ factors through one of the $\bigl[G\backslash
Z_{J,\delta}\bigr]$. We omit the details.

For each connected component $Z_{J,\delta}$ of~$Z'_J$ and for $\ubf \in
\Tscr(J)$ we have defined a $G$-invariant subscheme $Y^{\ubf}_J$ which
gives by passage to the quotient a locally closed substack
$\bigl[G\backslash Y^{\ubf}_J\bigr]$ of $\bigl[G\backslash
Z_{J,\delta}\bigr]$. The inverse images of these substacks $\Att^{\ubf}_0$
in~${\Att}_0$ are by definition the Ekedahl-Oort strata in ${\Att}_0$. Note
that these strata now carry a canonical scheme structure.

\secstart{} It is shown for $p > 2$ in~[\TWDimStrat] that $\pi$ is the
composition of a smooth morphism and a homeomorphism. In particular we obtain
that the codimension of the Ekedahl-Oort stratum $\Att^{\ubf}_0$
is the same as the codimension of $Y^{\ubf}$ in~$Z_{J,\delta}$ if
it is nonempty. Hence we get
$$
\codim(\Att^{\ubf}_0,\Att_0) = \dim(\Par_J) - \ell(u_\infty)
$$
by~\Dimensionformula. This gives a new proof of the main result
of~[\BMDFEO].

By~[\BMSTPEL]~3.2.7 the inverse image of the union of the open strata
in $\bigl[G\backslash Z'_J\bigr]$ is just the $\mu$-ordinary locus
of $\Att_0$ in the sense of [\TWOrdin] and we obtain a new proof of the
main result of~[\TWOrdin], as was already pointed out in~[\BMSTPEL].

\secstart{Example}:\label{\Fziptostronglat} {\it $F$-zips associated
to strongly divisible lattices.\/} Let $F$ be a finite extension
of~$\Qp$ with ring of integers~$O_F$ and finite residue field~$\kappa$
of cardinality~$q$. Let $k$ be a perfect extension of~$\kappa$ and let
$L = W(k) \otimes_{W(\kappa)} F$. Let $\sigma_0$ be the automorphism
of~$W(k)$ over~$W(\kappa)$ induced by the $q$th power Frobenius on~$k$,
and define $\sigma \in \Aut(L/F)$ by $\sigma = \sigma_0 \otimes \id_F$.
We fix a uniformizing element~$\pi$ of $O_F$, and we set $O_L := W(k)
\otimes _{W(\kappa)} O_F$.

Let $(H,\Phi,\Fil\updot)$ be a filtered isocrystal over~$L$. By this
we mean that $H$ is a finite dimensional $L$-vector space, equipped
with a $\sigma$-linear bijective operator $\Phi\colon H \to H$, and
with a descending filtration~$\Fil\updot$. Suppose that $\Mscr \subset
H$ is a strongly divisible lattice, i.e., an $O_L$-lattice such that
$$
\Mscr = \sum_{i\in\Z} \pi^{-i} \Phi\bigl(\Mscr \cap \Fil^i\bigr)\, .
\lformno\sublabel{\StrDivLatt}
$$
We claim that $M := \Mscr/\pi\Mscr$ naturally inherits the structure of
an $F$-zip over~$k$ with respect to the prime power~$q$. The definition
is as follows.

We let $C\updot$ be the descending filtration on~$M$ induced
by~$\Fil\udot$, so
$$
C^i = \bigl\{m \in M \bigm| \exists y \in \Mscr \cap \Fil^i\
\hbox{with}\ y \bmod \pi\Mscr = m \bigr\}\, .
$$
Next define the ascending filtration $D\dodot$ by
$$
D_i = \bigl\{m \in M \bigm| \exists y \in \Mscr\ \hbox{with}\
\pi^{-i}\Phi(y) \in \Mscr\ \hbox{and}\
\pi^{-i}\Phi(y) \bmod \pi\Mscr = m \bigr\}\, .
$$
Define a $k$-linear map $\tilde\phi_i\colon (C^i)^{(q)} \to D_i$ by
$\tilde\phi_i(m \otimes 1) = \pi^{-i}\Phi(y) \bmod \pi\Mscr$, where
$y \in \Mscr \cap \Fil^i$ is any element with $y \bmod \pi\Mscr = m$.
Using~\StrDivLatt\ we see that this is well-defined. It is easily seen
that $\tilde\phi_i$ vanishes on $(C^{i+1})^{(q)}$, so we may define
$\phi_i\colon (\gr^i_C)^{(q)} \to \gr_i^D$ to be the map induced
by~$\tilde\phi_i$.

It remains to be seen that $\phi_i$ is an isomorphism, and by a
dimension count it suffices to show that each~$\phi_i$ is surjective.
For this, consider an element $m \in D_i$. By assumption there is an
element $y \in \Mscr$ with $\pi^{-i}\Phi(y) \in \Mscr$ and
$\pi^{-i}\Phi(y) \bmod \pi\Mscr = m$. By~\StrDivLatt\ there are
elements $z_j \in \Mscr \cap \Fil^j$ such that $\pi^{-i} \Phi(y) =
\sum \pi^{-j} \Phi(z_j)$. Because $\Phi$ is injective, $\pi^{-i} y =
\sum \pi^{-j} z_j$. Let $y^\prime := y - \sum_{j<i} \pi^{i-j}z_j$.
Then $\Phi(y^\prime) \in \pi^i\Mscr$, and if we write $m^\prime$ for
the class of $\pi^{-i}\Phi(y^\prime)$ modulo~$\pi\Mscr$ then $m^\prime$
and~$m$ represent the same class in~$\gr_i^D$. On the other hand,
$y^\prime = \sum_{j\geq i} \pi^{i-j} z_j$ lies in $\Fil^i \cap \Mscr$.
Hence it follows that $\bar{m} \in \gr_i^D$ is in the image of~$\phi_i$.

If we look more closely, we see that we not only get the structure of
an $F$-zip on~$M$, but that it comes equipped with a natural splitting
of the $D\dodot$-filtration. Namely, define $\tilde\Mscr \subset H$ by
$\tilde\Mscr := \sum_{i\in\Z} \pi^{-i}(\Mscr \cap \Fil^i)$. Then $\Phi$
gives an isomorphism $\sigma^* \tilde\Mscr \arriso \Mscr$. Moreover, we
have a natural isomorphism $\tilde\Mscr/\pi \tilde\Mscr \arriso
\oplus_{i\in\Z} \gr^i_C(M)$. Hence $\Phi$ induces an isomorphism $\oplus
(\gr^i_C)^{(q)} \arriso M$; it maps $(\gr^i_C)^{(q)}$ into~$D_i$, and
the composition $(\gr^i_C)^{(q)} \to D_i \to \gr_i^D$ is of
course~$\phi_i$. The image of $(\gr^i_C)^{(q)}$ inside~$D_i$ is the
subspace $E_i \subset D_i$ given by $E_i = \bigl\{m \in M \bigm| \exists
y \in \Mscr \cap \Fil^i\ \hbox{with}\ \pi^{-i}\Phi(y) \bmod \pi\Mscr =
m \bigr\}$, and this is a complement for~$D_{i-1}$ inside~$D_i$.

\secstart{Example}: {\it $F$-zips associated to K3 surfaces.\/}
Fix a natural number~$d$ and a prime number~$p$ with $p \nmid 2d$.
Let $S$ be a scheme of characteristic~$p$, and let $(Y,L)$ be a 
K3 surface with polarization of degree~$2d$ over~$S$. This gives 
rise to a sequence
$$
\matrix{ S = S_1 \supset S_2 \supset \cdots \supset S_{10} \supset
\kern-1.2em & S_{\infty} &\cr &\parallel&\cr
&S_{\infty,10}& \kern-.7em\supset S_{\infty,9} \supset \cdots
\supset S_{\infty,1}\, .\cr}\lformno\sublabel{\KThreeStrat}
$$
Here $S_h \subset S$, for $h \in \{1,2,\ldots,10,\infty\}$, is the
closed subscheme of~$S$ given, loosely speaking, by the condition
that the formal group of~$X$ has height $\geq h$. For details we
refer to Ogus's paper~[\OgusKKK]. On the supersingular locus 
$S_\infty$ we have a further set-theoretic stratification, letting
$S_{\infty,i}$ be the locus of points $s \in S$ where $\sigma_0(Y_s) 
\leq i$; here $\sigma_0$ denotes the Artin invariant. For details 
we again refer to~[\OgusKKK].

We now want to connect the stratification in~\KThreeStrat\ with
our theory of $F$-zips. Let $f\colon Y \to S$ be the structural
morphism, and consider the second de Rham cohomology $H :=
R^2f_\ast \Omega\udot_{Y/S}$. It is a locally free
$\Oscr_S$-module of rank~$22$, which comes equipped with a
non-degenerate symmetric bilinear form $Q\colon H \times H \to
\Oscr_S$. If $c_1(L) \in H(S)$ is the first Chern class of~$L$
then the primitive cohomology
$$
M := \bigl\langle c_1(L) \bigr\rangle^\perp \subset H
$$
is locally free of rank~$21$, and $Q$ restricts to a non-degenerate
form on~$M$, which we again call~$Q$.

As in~\HodgeWitttoFzip, let $C\udot$ be the Hodge filtration
on~$M$ and let $D\dodot$ be the conjugate filtration (up to a
renumbering). The type~$\tau$ is given by $\tau(0)=\tau(2) = 1$
and $\tau(1) = 19$, and $C\udot$ and~$D\dodot$ are orthogonal
filtrations. The inverse Cartier isomorphism gives isomorphisms
$\phi_i \colon (\gr^i_C)^{(p)} \arriso \gr_i^D$ such that $\Mline
= (M,Q,C\udot,D\dodot,\phi\dodot)$ is an orthogonal $F$-zip.

Let $(V,\psi)$ be an orthogonal space over~$\Fp$ with $\dim(V) =
21$. Set $G := \SO(V,\psi)$, which has root system of
type~B$_{10}$. We take a basis of simple roots as
$\{\alpha_1,\ldots,\alpha_{10}\}$ as in [\Bourb], Planche~II; thus,
$\alpha_{10}$ is the short root. Let $I = \{s_1,\ldots,s_{10}\}$
be the corresponding set of simple reflections. We have
$$
W \cong \bigl\{\rho \in S_{21} \bigm| \rho(j) + \rho(22-j) = 22
\quad \hbox{{\rm for all $j$}}\bigr\}\, ,
$$
with $s_i$ corresponding to the element $(i,i+1) \cdot
(21-i,22-i)$ for $1 \leq i \leq 9$ and $s_{10}$ corresponding to
the transposition $(10,12)$. (Note that $\rho(11) = 11$ for all
$\rho \in W$.)

Let $J := I\setminus\{s_1\}$. We have $W_J = \bigl\{\rho \in
W\bigm| \rho(1) = 1\bigr\}$, so $W_J\backslash W$ is a set of $20$
elements. The set ${}^JW$ of minimal representatives consists of
elements $x_1,\ldots,x_{20}$, which we number in such a way that
$\ell(x_j) = 20-j$. In the Bruhat ordering we have $x_1 > x_2 >
\cdots > x_{20}$.

Consider the covering ${}^\#S \to S$, such that the $T$-valued 
points of~${}^\#S$, for $T$ an $S$-scheme, are the isometries 
$M_T \arriso V_T$. Then ${}^\#S$ is an torsor over~$S$ in the 
\'etale topology, under the group $\OO(V,\psi)$. If $X_\tau$ is 
the scheme of orthogonal $F$-zip structures on $(V,\psi)$, as 
in~\SymplOrth, then we have an $\OO(V,\psi)$-equivariant morphism 
$\rho\colon {}^\#S \to X_\tau$. If $X_\tau^{(j)} \subset X_\tau$ 
is the stratum corresponding to the element~$x_j$, let $S^{(j)} 
\subset S$ be the subscheme obtained as the quotient of 
$\rho^{-1}\bigl(X_\tau^{(j)}\bigr)$ under $\OO(V,\psi)$.

The connection between the EO-loci $S^{(j)}$ thus obtained and the strata
in~\KThreeStrat\ is given by the following result.

\secstart{Proposition}:\label{\KKKStratifs} {\sl With notation 
as in~\KThreeStrat\ we have
$$
S^{(j)} = S_j \setminus S_{j+1}\qquad \hbox{for $1 \leq j \leq 11$}\, ,
$$
where we let $S_{11} := S_\infty = S_{\infty,10}$ and $S_{12} :=
S_{\infty,9}$. Further, on $k$-valued points we have
$$
S_{\infty,l}(k) = \coprod_{i \geq 21-l} S^{(i)}(k) \qquad
\hbox{for $1 \leq l \leq 10$}\, .
$$}

This is essentially what Ogus proves in~[\OgusKKK]. Thus, we see
that our theory of $F$-zips gives a uniform and scheme-theoretic
approach to the whole chain in~\KThreeStrat.


\vskip 10mm

\noindent
{\partifont References}
\secskip

\item{[\BorTits]} A.~Borel, J.~Tits: {\it Groupes r\'eductifs\/}, 
Publ.\ Math.\ de l'IHES {\bf 27} (1965), 55--151.
\item{[\Bourb]} N.~Bourbaki: {\it Groupes et alg\`ebres de Lie, 
Chap.\ 4,5 et 6\/}, Masson, Paris, 1981.
\item{[\DelCohInt]} P.~Deligne: {\it Cohomologie des intersections 
compl\`etes\/}, SGA7, Exp.~XI, LNM {\bf 340}, Springer-Verlag, 
Berlin, 1973, pp.~39--61.
\item{[\DelVarSh]} P.~Deligne: {\it Vari\'et\'es de Shimura: 
interpr\'etation modulaire, et techniques de construction de 
mod\`eles canoniques\/}, in: A.~Borel and W.~Casselman (eds.), 
Automorphic Forms, Representations, and $L$-functions (2),  Proc.\ 
of Symp.\ in Pure Math. {\bf 33}, AMS, Providence, RI, 1979, 
pp.~247--290.
\item{[\DelRelKKK]} P.~Deligne: {\it Rel\`evement des surfaces K3 en 
caract\'eristique nulle (r\'edige par L.~Illusie)\/}, in: 
Surfaces Alg\'ebriques, J.~Giraud et al., eds., LNM {\bf 868}, 
Springer-Verlag, Berlin, 1981, pp.~58--79.
\item{[\DelIll]} P.~Deligne, L.~Illusie: {\it Rel\`evements modulo 
$p^2$ et d\'ecomposition du complexe de de Rham\/}, Invent.\ 
Math.\ {\bf 89} (1987), 247--270.
\item{[\SGAthree]}
M.~Demazure et al.: {\it Sch\'emas en groupes, I,
II, III\/}, LNM {\bf 151}, {\bf 152}, {\bf 153}, 
Springer-Verlag, Berlin, 1970.
\item{[\MAV]} C.~Faber, G.\ van der Geer, F.~Oort (eds.): {\it 
Moduli of abelian varieties\/}, Progr.\ in Math.\ {\bf 195}, 
Birkh\"auser Verlag, Basel, 2001.
\item{[\GeerKats]} G.~van der Geer, T.~Katsura: {\it An invariant
for varieties in positive characteristic\/}, in:  S.~Vostokov and 
Yu.~Zarhin (eds.), Algebraic number theory and algebraic geometry, 
Contemp.\ Math.\ {\bf 300}, AMS, Providence, RI, 2002, pp.~131--141. 
\item{[\EGA]} A.~Grothendieck: {\it \'El\'ements de g\'eom\'etrie
alg\'ebrique (r\'edig\'es avec la collaboration de 
J.~Dieu\-donn\'e)\/}, Publ.\ Math.\ de l'IHES {\bf 4}, {\bf 8}, 
{\bf 11}, {\bf 17}, {\bf 20}, {\bf 24}, {\bf 28}, {\bf 32} 
(1960--67).
\item{[\IllRay]} L.~Illusie, M.~Raynaud: {\it Les suites spectrales
associ\'ees au complexe de de Rham-Witt\/}, Publ.\ Math.\ de 
l'IHES {\bf 57} (1983), 73--212.
\item{[\KatoLog]} K.~Kato: {\it Logarithmic structures of 
Fontaine-Illusie\/}, in: J.-I.~Igusa (ed.), Algebraic analysis, 
geometry, and number theory, Johns Hopkins Univ.\ Press, 
Baltimore, MD, 1989, pp.~191--224
\item{[\KatzDiffEq]} N.~Katz: {\it Algebraic solutions of differential 
equations ($p$-curvature and the Hodge filtration)\/}, Invent.\ 
Math.\ {\bf 18} (1972), 1--118.
\item{[\KottPts]} R.~Kottwitz: {\it Points on some Shimura varieties 
over finite fields\/}, J.A.M.S.\ {\bf 5} (1992), 373--444.
\item{[\Kraft]} H.~Kraft: {\it Kommutative algebraische $p$-Gruppen 
(mit Anwendungen auf $p$-divisible Gruppen und abelsche Variet\"aten)\/}, 
manuscript, Univ.\ Bonn, Sept.\ 1975, 86 pp. (Unpublished)
\item{[\LusztPar]} G.~Lusztig: {\it Parabolic character sheaves I, II\/}, 
preprints 2003, math.RT/0302151 and\break math.RT/0302317.
\item{[\Mats]} H.~Matsumura: {\it Commutative ring theory\/}, 
Cambridge Studies in Adv.\ Math.\ {\bf 8}, Cambridge Univ. Press,
Cambridge, 1986.
\item{[\BMGSAS]} B.~Moonen: {\it Group schemes with additional 
structures and
Weyl group cosets\/}, in: [\MAV], pp.~255--298.
\item{[\BMSTPEL]} B.~Moonen: {\it Serre-Tate theory for moduli spaces 
of
PEL type\/}, to appear in the Ann.\ scient.\ \'Ec.\ Norm.\ 
Sup.
\item{[\BMDFEO]} B.~Moonen: {\it A dimension formula for Ekedahl-Oort
strata\/}, preprint 2002, math.AG/0208161, to appear in
the Ann.\ Inst.\ Fourier (Grenoble).
\item{[\Oda]} T.~Oda: {\it The first de Rham cohomology group and 
Dieudonn\'e modules\/}, Ann.\ scient.\ \'Ec.\ Norm.\ Sup.\ (4), 
{\bf 2}, (1969), 63--135.
\item{[\OgusKKK]} A.~Ogus: {\it Singularities of the height strata in 
the moduli of K3 surfaces\/}, in: [\MAV], pp.~325--343.
\item{[\OortStrat]} F.~Oort: {\it A stratification of a moduli space of 
polarized abelian varieties in positive characteristic\/}, in: 
C.~Faber and E.~Looijenga (eds.), Moduli of curves and abelian 
varieties, Aspects Math.\ {\bf E33}, Vieweg, Braunschweig, 1999, 
pp.~47--64.
\item{[\OortTexel]} F.~Oort: {\it A stratification of a moduli space of 
abelian varieties\/}, in: [\MAV], pp.~435--416.
\item{[\RapZink]} M.~Rapoport, T.~Zink: {\it Period spaces for 
$p$-divisible groups\/},
Annals of Math.\ Studies {\bf 141}, 
Princeton Univ.\ Press, Princeton, NJ, 1996.
\item{[\TWOrdin]} T.~Wedhorn: {\it Ordinariness in good reductions of 
Shimura varieties of PEL-type\/}, Ann.\ scient.\ \'Ec.\ Norm.\ 
Sup.\ (4), {\bf 32}, (1999), 575--618.
\item{[\TWDimStrat]} T.~Wedhorn: {\it The dimension of Oort strata of 
Shimura varieties of PEL-type\/}, in: [\MAV], pp.~441--471.

\bigskip

Ben Moonen, University of Amsterdam, KdV Institute for Mathematics,
Plantage Muidergracht 24, 1018 TV, Amsterdam, The Netherlands.
E-mail: {\tt bmoonen@science.uva.nl}

\smallskip

Torsten Wedhorn, University of Bonn, Mathematisches Institut Bonn,
Beringstra\ss e 4, 53115 Bonn, Germany.
E-mail: {\tt wedhorn@math.uni-bonn.de}

\bye